\documentclass{article}
\usepackage[utf8]{inputenc}

\usepackage{amssymb}
\usepackage{amsmath,mathtools}
\usepackage{amsthm}
\usepackage{amsfonts}
\usepackage{tikz}
\usepackage{xcolor}
\usepackage[top=1.0in, bottom=1.0in, left=1.0in, right=1.0in, headheight=1.0in]{geometry}
\usepackage{mathtools}
\usepackage{youngtab}
\usepackage{graphicx}
\usepackage{multicol}
\usepackage{paracol}
\graphicspath{ {./Graphs/} }
\usepackage[nottoc,notlof,notlot]{tocbibind}
\usepackage{changepage}

\newtheorem{theorem}{Theorem}[section]
\newtheorem{lemma}[theorem]{Lemma}

\newtheorem{prop}[theorem]{Proposition}

\newtheorem{corollary}[theorem]{Corollary}
\newtheorem{remark}[theorem]{Remark}

\newcommand{\hng}{$1$-$HNG$}
\newcommand{\ahng}{$a$-$HNG$}
\newcommand{\zerong}{$0$-$NG$}
\newcommand{\oneng}{$1$-$NG$}
\newcommand{\twohng}{$2$-$HNG$}
\newcommand{\zerohng}{$0$-$HNG$}
\newcommand{\ang}{$a$-$NG$}
\newcommand{\hngsp}{$1$-$HNG$ }
\newcommand{\ahngsp}{$a$-$HNG$ }
\newcommand{\zerongsp}{$0$-$NG$ }

\newcommand{\zerohngsp}{$0$-$HNG$ }
\newcommand{\angsp}{$a$-$NG$ }

\newcommand{\aplusonehng}{$(a+1)$-$HNG$}

\newcommand{\aplusoneng}{$(a+1)$-$NG$}

\date{\today  }
\title{Hereditary Nordhaus-Gaddum Graphs}

\author{Vaidy Sivaraman\\
\small Dept. of Mathematics and Statistics\\
\small Mississippi State University\\
\small Starkville, MS 39762\\
\small\tt vs692@msstate.edu\
\and
Rebecca Whitman\\
\small Department of Mathematics\\
\small University of California, Berkeley\\
\small Berkeley, CA 94720\\
\small\tt  rebecca\_whitman@berkeley.edu\
}
\begin{document}
\maketitle

\begin{abstract}    
Nordhaus and Gaddum proved in 1956 that the sum of the chromatic number $\chi$ of a graph $G$ and its complement is at most $|G|+1$. The Nordhaus-Gaddum graphs are the class of graphs satisfying this inequality with equality, and are well-understood. In this paper we consider a hereditary generalization: graphs $G$ for which all induced subgraphs $H$ of $G$ satisfy $\chi(H) + \chi(\overline{H}) \le |H|$. We characterize the forbidden induced subgraphs of this class and find its intersection with a number of common classes, including line graphs. We also discuss $\chi$-boundedness and algorithmic results. 
  \end{abstract}
  
  
  \smallskip
  

\section{Introduction}

In this paper we define a new family of graph classes generalizing the Nordhaus-Gaddum graphs, and introduce the largest hereditary subclasses of each of these classes. Throughout, all graphs are finite and simple, with neither loops nor multiple edges. We begin by introducing four graph invariants used throughout the paper. Let $G = (V(G), E(G))$ be a graph. Let $n$ denote $|V(G)|$ and $\overline{G}$ denote the complement of $G$. The \textit{clique number} $\omega(G)$ is the largest positive integer $k$ such that the complete graph on $n$ vertices, denoted $K_n$, is a subgraph of $G$, and the \textit{independence number} $\alpha(G)$ is the largest positive integer $l$ such that $\overline{K_n}$ is a subgraph of $G$. The \textit{chromatic number} $\chi(G)$ is the smallest positive integer $k$ such that the vertices of $G$ can be properly colored with $k$ colors (a proper coloring is one where adjacent vertices receive different colors). The \textit{clique cover number} $\theta(G)$ is the smallest positive integer $l$ such that $V(G)$ can be partitioned into $l$ sets, each inducing a clique. It follows from the definitions that $\omega(\overline{G}) = \alpha(G)$, $\chi(\overline{G}) = \theta(G)$, $\omega(G) \le \chi(G) \le n$, and $\alpha(G) \le \theta(G) \le n$. 

The relationship between the chromatic number of a graph and that of its complement was first studied by Nordhaus and Gaddum in 1956 \cite{NoGa56}. They proved that the sum of the two invariants is bounded above by $n + 1$: 
\begin{equation}
\label{eq:NG_inequality}
    \chi(G) + \theta(G) \le n + 1
\end{equation}

Since then several researchers have studied the same problem for various graph invariants. See \cite{AoHa13} for a survey. The class of graphs satisfying Inequality \ref{eq:NG_inequality} with equality are called the Nordhaus-Gaddum graphs (we denote this class \zerong). The class \zerongsp has been widely studied, with three increasingly elegant structural characterizations \cite{Fi66} \cite {StTu08} \cite{CoTr13}. Cheng, Collins, and Trenk gave a degree sequence characterization of \zerongsp and an enumeration of the class \cite{ChCoTr16}. 

Notably, \zerongsp is not a hereditary class, closed under taking induced subgraphs. As an example, the cycle graph on $5$ vertices, $C_5$, has $\chi = \theta = 3$ and $n = 5$, so satisfies Inequality \ref{eq:NG_inequality} with equality. Its induced subgraph the path graph $P_4$ has $\chi = \theta = 2$ and $n = 4$, so is not in \zerong. A natural question is to understand the largest hereditary subclass of \zerong, i.e. the set of graphs for which all induced subgraphs are Nordhaus-Gaddum. It turns out that this class coincides with a class well-studied in the context of induced subgraphs, namely, the threshold graphs. We provide a proof in Proposition \ref{prop:threshold}, later in this section. Threshold graphs have been relaxed and generalized in several different ways (see for instance \cite{MaPe95} \cite{YCC} \cite{BSZ}), and this paper will give another. 

The class \zerongsp is small and contains few important subclasses, aside from the threshold graphs. Nevertheless, a number of classes have interesting intersections. One is the $C_4, 2K_2$-free graphs (the pseudo-split graphs), defined in \cite{Bl93}. Pseudo-split graphs containing an induced $C_5$ satisfy Inequality \ref{eq:NG_inequality} with equality, and form one of three subclasses of \zerongsp in Collins and Trenk's characterization \cite{CoTr13}. Bl\'{a}szik \cite{Bl93} showed all pseudo-split graphs satisfy the relaxed bound $\chi(G) + \theta(G) \ge n$. 

We build on this generalization to introduce a parameter $a$ to Inequality \ref{eq:NG_inequality} and consider graphs satisfying the following inequality for each choice of $a \ge 0$: 
\begin{equation}
\label{eq:generalized_ng}
    \chi(G) + \theta(G) \ge n + 1 - a. 
\end{equation}

Given $a \ge 0$, a graph $G$ is \textit{a-Nordhaus-Gaddum} if $G$ satisfies Inequality \ref{eq:generalized_ng} with equality. The generalized classes of $a$-Nordhaus-Gaddum graphs are again not hereditary; as an example, $C_{2a+5}$ has $\chi = 3$ and $\theta = a+3$, so is an element of \ang, but its induced subgraph $P_{2a+4}$ has $\chi = 2$ and $\theta = a+2$, so is not an element of \ang. We hence define $G$ to be \textit{a-hereditary-Nordhaus-Gaddum} if every induced subgraph $H$ of $G$ satisfies \ref{eq:generalized_ng} with equality. We denote the class of $a$-Nordhaus-Gaddum graphs by \angsp and the class of $a$-hereditary-Nordhaus-Gaddum graphs by \ahng. For all $a \ge 0$, \angsp and \ahngsp are both closed under complementation, and \ahngsp is hereditary. 

We have the following chain of inclusions: 

\begin{prop}
\label{prop:inclusions}
    $\text{\zerohng} \subset \text{\zerong} \subset \text{\hng} \subset \text{\oneng} \subset \text{\twohng} \subset \ldots$
\end{prop}
\begin{proof}
    Fix $a \ge 0$. By definition \ahngsp $\subset$ \ang. Suppose for a contradiction that \angsp $\not \subset $ \aplusonehng, so there exists $G \in $ \angsp with an induced subgraph $H \not \in$ \aplusoneng. Let $K$ be the induced subgraph with vertex set $V(G) - V(H)$. Then $\chi(G) + \theta(G) \ge |G| + 1 - a$ and $\chi(H) + \theta(H) < |H| + 1 - (a+1)$. By Inequality \ref{eq:NG_inequality}, $\chi(K) + \theta(K) \le |K| + 1$. As a result: 
    $$\chi(G) + \theta(G) \le \chi(H) + \chi(K) + \theta(H) + \theta(K) < |H| + |K| + 1 - a = |G| + 1 - a, $$
    contradicting $G$'s inclusion in \ang. Thus instead \angsp $\subset$ \aplusonehng. 

    Since $C_{2a+5}$ is in \angsp but not \ahng, and $P_{2a+4}$ is in \aplusonehng but not \ang, all inclusions are strict. 
\end{proof}

We introduce more notation used throughout the paper, then discuss a number of important graph classes and their relation to \angsp and \ahng. Given positive integers $k, l$, we denote the complete bipartite graph with $k$ vertices in one half of the bipartition and $l$ in the other by $K_{k,l}$. We use $+$ to denote the disjoint union of two graphs. Given a vertex $v$ (resp. induced subgraph $H$), let $G - \{v\}$ (resp., $G - H$) denote the induced subgraph on vertex set $V(G) - \{v\}$ (resp., $V(G) - V(H)$).  A vertex $v$ is \textit{complete to} $H \subseteq V(G)$ if $v$ is adjacent to all vertices in $H$. A vertex subset $A$ is complete to $H$ if all $v \in A$ are complete to $H$. The neighborhood of a vertex $v$ in $G$ (resp., in induced subgraph $H$) is denoted $N(v)$ (resp., $N_H(v)$). A vertex $v$ is \textit{isolated} in $G$ (resp., in $H$) if $N(v)$ (resp. $N_H(v)$) is empty, and \textit{dominating} if $N(v) = V(G) - \{v\}$. We refer to two induced subgraphs as isolated from one another where there are no edges with one endpoint in each subgraph. 

A graph $G$ is said to be \emph{threshold} if there exist a function $w: V(G) \to \mathbb{R}$ and a real number $t$ such that there is an edge between two distinct vertices $u$ and $v$ if and only if $w(u) + w(v) > t$ \cite{ChHa73}. The class of threshold graphs has been studied in great detail. See Mahadev and Peled's book for more information \cite{MaPe95}. Threshold graphs are identifiable by forbidden induced subgraphs and with a vertex ordering (\cite{ChHa77}), among myriad other characterizations. 

\begin{theorem} 
\label{prop:threshold_char}
Given a graph $G$, the following are equivalent: 
\begin{enumerate}
    \item[(i)] $G$ is threshold. 
    \item[(ii)] $G$ contains no induced $2K_2, P_4, C_4$. 
    \item[(iii)] $V(G)$ can be given an ordering $v_1, \ldots, v_n$ such that for every $i$, $1 \le i \le n$, $v_i$ is adjacent to all or none of $v_1, \ldots, v_{i-1}$. 
\end{enumerate} 
\end{theorem}

Together, these two characterizations show that threshold graphs are exactly \zerohng. 

\begin{prop}
\label{prop:threshold}
    The class \zerohngsp is exactly the class of threshold graphs. 
\end{prop}
\begin{proof}
    Given $a \ge 0$ and a graph $G$ with an isolated or dominating vertex $v$, $G$ satisfies Inequality \ref{eq:generalized_ng} with equality if and only if $G-\{v\}$ satisfies Inequality \ref{eq:generalized_ng} with equality. Hence Theorem \ref{prop:threshold_char} (iii) implies that threshold graphs are contained in \zerohng. For the converse it suffices to see from Theorem \ref{prop:threshold_char}(ii), that none of $P_4, 2K_2, C_4$, the three forbidden induced subgraphs for the class of threshold graphs, is in \zerohng. 
\end{proof}
 
A graph $G$ is \textit{split} if its vertices can be partitioned into a clique and a stable set \cite{FoHa77}. A graph $G$ is \emph{chordal} if every induced cycle in it is a triangle \cite{Di61}. It is \emph{weakly chordal} if every induced cycle in it or its complement is either a triangle or a square \cite{Ha85}. A graph $G$ is \emph{perfect} if $\chi(H) = \omega(H)$ for all induced subgraphs $H$ of $G$ \cite{Be61}. The Strong Perfect Graph Theorem \cite{CRST} shows a graph is perfect if and only if it contains none of $\{C_5, C_7, \overline{C_7}, C_9,  \overline{C_9}, \ldots\}$ as an induced subgraph. 

We also define a new class generalizing perfect graphs: a graph is \textit{apex-perfect} if it contains a vertex whose deletion results in a perfect graph. Although hereditary, the forbidden induced subgraph characterization of this class remains unknown. 

By the forbidden induced subgraph characterizations, we have the following chain of inclusions on graph classes: 
    $$\text{Threshold} \subset \text{Split} \subset \text{Chordal} \subset \text{Weakly Chordal} \subset \text{Perfect} \subset \text{Apex-Perfect}.$$

We compare the hereditary classes \ahngsp to the above hereditary classes. First, \zerohngsp is exactly the class of threshold graphs. In a split graph, $\omega(G) + \alpha(G) \ge |G|$ \cite{HaSi81}, so since split graphs are perfect, $\chi(G) + \theta(G) = \omega(G) + \alpha(G) \ge |G|$, it follows that split graphs are a subclass of \hng. For $a \ge 1$, \ahngsp is not a subclass of chordal, weakly chordal, or perfect graphs, since it contains $C_5$. The converse also holds: chordal, weakly chordal, and perfect graphs do not constitute subclasses of \ahng, since $P_{2a+4}$ is in each of these classes but not in \ahng. We prove later in Theorem \ref{prop:apex_perfect} that all graphs in \hngsp are apex-perfect. 

Another generalization of perfect graphs is to $\chi$-bounding functions. For a hereditary class $\mathcal{G}$, $f: \mathbb{N} \rightarrow \mathbb{N}$ is a $\chi$-bounding function on $\mathcal{G}$ if for all $G \in \mathcal{G}$, $\chi(G) \le f(\omega(G))$. For perfect graphs, the function $f(x) = x$ is $\chi$-bounding, and this is of course best possible. See \cite{ScSe20} for a survey of hereditary classes with known $\chi$-bounding functions. Since all graphs in \hngsp are apex-perfect, it follows as Theorem \ref{prop:chi_bound_Vizing} that \hngsp is $\chi$-bounded by the function $f(x) = x+1$, which is best possible. 

Given a graph $G$, it holds that $\omega(G) + \alpha(G) \le n + 1$. Where this bound holds with equality for all induced subgraphs $H$ of $G$, it follows that $G$ is threshold \cite{LiPoSi19}. More interesting is the generalization of threshold graphs to the class of sum-perfect graphs. A graph is \textit{sum-perfect} if for all induced subgraphs $H$, $\omega(H) + \alpha(H) \ge n$. One of the authors, together with Litjens and Polak, provides a forbidden induced subgraph characterization of the class \cite{LiPoSi19}. It follows from the characterization that the class of sum-perfect graphs is a subclass of the weakly-chordal graphs, and hence perfect. Since $\omega(G) \le \chi(G)$ and $\alpha(G) \le \theta(G)$, sum-perfect graphs are also \hng. We make use of this inclusion in Theorem \ref{prop:FIS_for_1_hng}, our forbidden induced subgraph characterization of \hng. 

The line graph of $G$, denoted $L(G)$, has the edges of $G$ as its vertices. Two vertices in $L(G)$ are adjacent if as edges in $G$ they are incident to a shared vertex. Line graphs translate questions about edges into questions about vertices, and so are an important tool for simplifying otherwise-intractable problems. The class of line graphs is hereditary and there is a beautiful characterization of this class in terms of its forbidden induced subgraphs \cite{LWB}. We characterize the intersection of \hngsp and line graphs in Theorem \ref{prop:hng_line_graph}. 

A \emph{claw} is an induced subgraph that is isomorphic to $K_{1,3}$; its trivalent vertex is called its \emph{center}. A graph is \emph{claw-free} if it contains no claw. Since the claw is one of the forbidden induced subgraphs of line graphs, claw-free graphs are an important generalization of line graphs. They also have strong algorithmic properties (see the survey \cite{CFGS}). We characterize the intersection of \hngsp and claw-free graphs in Theorem \ref{prop:hng_claw_free}. Section 5 ends with characterizations of the intersection of \hngsp with two simple, essential classes:  bipartite graphs and triangle-free graphs. 

The paper is structured as follows. In Section 2 we study the effect of vertex deletion on $\chi$ and $\theta$, which is connected to forbidden induced subgraphs of \ahng. From here, the remaining sections of the paper pertain only to the smallest class \hng. In Section 3 we extend the forbidden induced subgraph characterization of sum-prefect graphs \cite{LiPoSi19} to characterize the forbidden induced subgraphs of \hng. In Section 4 we show graphs in \hngsp are apex-perfect, thereby giving a best possible $\chi$-bounding function on the class. Section 5 centers on the intersection of \hngsp with, respectively, the classes of line graphs, claw-free graphs, and triangle-free graphs. In Section 6 we provide a number of optimization results about \hng, showing that inclusion in \hng, $\omega(G), \alpha(G), \chi(G),$ and $\theta(G)$ can all be computed in polynomial time. 

\section{Vertex Deletion in Minimum Colorings and Clique Coverings}

In this section we present several results about vertex deletion and forbidden induced subgraphs of \ahng. First, since \ahngsp is closed under complementation, so too is its set of forbidden induced subgraphs, which we record as a remark. 

\begin{remark}
\label{prop:FIS_closed_complement}
    For all $a \ge 0$, the set of forbidden induced subgraphs of \ahngsp is closed under complementation. 
\end{remark}

Second, each vertex in a graph can contribute at most one color to a coloring and add at most one clique to a clique covering. Since a proper coloring of a graph $G$ is formally a partition of $V(G)$ into stable sets, we call each stable set a \textit{color class}, that is, the set of vertices with the same color. A vertex $v$ is \textit{$\chi$-distinct} in $G$ if there exists a minimum proper coloring of $G$ where $v$ is the only vertex in its color class. Similarly, $v$ is \textit{$\theta$-distinct} in $G$ if there exists a minimum proper clique covering where $v$ is the only vertex in its clique. Such a minimum proper coloring or clique covering is called \textit{$v$-distinct}. Deleting a vertex $v$ from $G$ will decrement $\chi$ or $\theta$ if and only if $v$ is $\chi$-distinct or $\theta$-distinct, respectively, which we record in the following proposition. 

\begin{prop}
\label{prop:vertex_deletion}
    Given $v \in G$, $\chi(G - \{v\}) = \chi(G) - 1$ if and only if $v$ is $\chi$-distinct. Otherwise $\chi(G - \{v\}) = \chi(G)$. Analogously, $\theta(G - \{v\}) = \theta(G) - 1$ if and only if $v$ is $\theta$-distinct, and otherwise $\theta(G - \{v\}) = \theta(G)$. 
\end{prop}
\begin{proof}
    Let $v \in G$. Clearly $\chi(G - \{v\}) = \chi(G)$ or $\chi(G) - 1$. If $v$ is $\chi$-distinct, then choose a $v$-distinct coloring. The restriction of this coloring to $G - \{v\}$ gives a proper coloring of $G - \{v\}$ with $\chi(G) - 1$ colors, so $\chi(G - \{v\}) \le \chi(G) - 1$. Hence they are equal. For the reverse direction, if $\chi(G - \{v\}) = \chi(G) - 1$, then fix a proper coloring of $G - \{v\}$ using $\chi(G) - 1$ colors. Add a new color class containing only $v$ to produce a proper coloring of $G$ using $\chi(G)$ colors. Hence the coloring is minimal in $G$, and we conclude that $v$ is $\chi$-distinct. An analogous proof holds for $\theta(G)$. 
\end{proof}

The following result is used to identify $\theta$-distinct vertices; an analogous result holds for $\chi$-distinct vertices. 

\begin{prop}
\label{prop:identifying_deletions}
    For any $v \in V(G)$ such that $N(v)$ induces a stable set, $v$ is $\theta$-distinct in $G$ if its neighbors are not $\theta$-distinct in $G - \{v\}$. 
\end{prop}
\begin{proof}
    Let $v \in V(G)$ with $N(v)$ inducing a stable set. We prove the contrapositive. Suppose $v$ is not $\theta$-distinct in $G$. Then for all minimum proper clique coverings of $G$, $v$ shares a clique with one of its neighbors. Choose a minimum proper clique covering. This clique covering restricted to $G - \{v\}$ is then a proper clique covering with some $w \in N(v)$ in a distinct clique. Proposition \ref{prop:vertex_deletion} implies that $\theta(G - \{v\}) = \theta(G)$, so the clique covering restricted to $G - \{v\}$ is minimum. Thus $w$ is $\theta$-distinct in $G - \{v\}$. 
\end{proof}

Identifying distinct vertices is critical because forbidden induced subgraphs cannot contain distinct vertices of either type. 

\begin{prop}
\label{prop:vertex_deletion_mfis}
    If $G$ is a forbidden induced subgraph of \ahng, then $G \in$ \aplusonehng, and no vertex of $G$ is $\chi$-distinct or $\theta$-distinct.
\end{prop}
\begin{proof}
    Let $v \in V(G)$. Since $G$ is a forbidden induced subgraph of \ahng, it follows that $G - \{v\} \in$ \ahng. Thus $\chi(G) + \theta(G) \le n-a$ and $\chi(G - \{v\}) + \theta(G - \{v\}) \ge n-a$. Since $\chi(G) \ge \chi(G - \{v\})$ and $\theta(G) \ge \theta(G - \{v\})$, it follows that $\chi(G) = \chi(G-\{v\})$, $\theta(G) = \theta(G - \{v\}$, and $\chi(G) + \theta(G) = n-a$. Thus $G \in $ \aplusonehng, and by Proposition \ref{prop:vertex_deletion}, $v$ is not $\chi$-distinct or $\theta$-distinct. 
\end{proof}

The above two propositions combine to give the following corollary, which we use extensively in our characterization of \hngsp in Section 3. 

\begin{corollary}
\label{prop:theta_distinct_mfis}
    If there exists $v \in V(G)$ such that $N(v)$ induces a stable set and no vertex in $N(v)$ is $\theta$-distinct in $G - \{v\}$, then $G$ is not a forbidden induced subgraph of \ahng for any $a \ge 0$. 
\end{corollary}

\section{A Forbidden Induced Subgraph Characterization of \hng}

\begin{figure}[ht]
\centering
  \includegraphics[height=10cm]{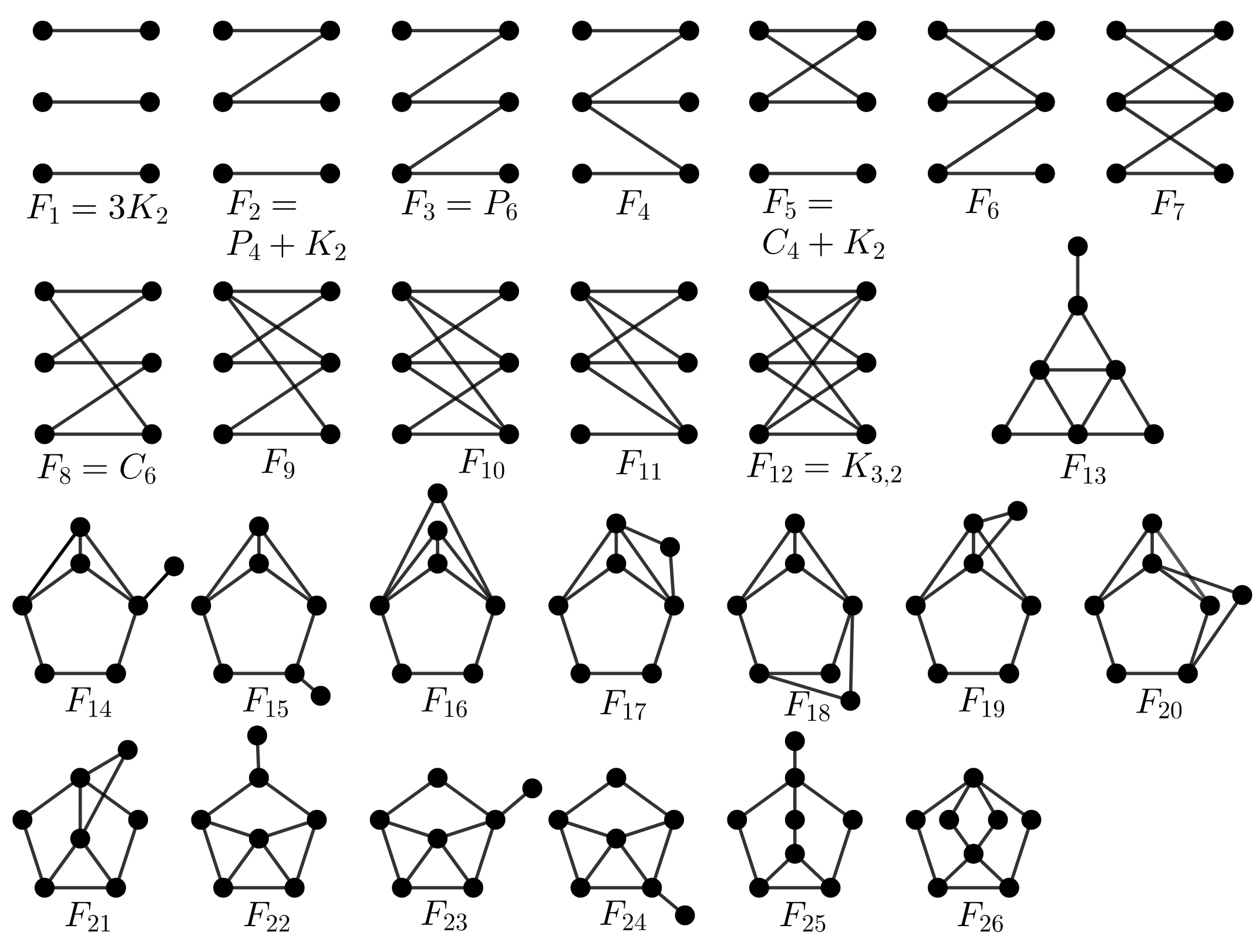}
  \caption{These graphs, together with their complements, comprise the set $\mathcal{F}$ of forbidden induced subgraphs of \hng.}
  \label{fig:fis_for_1hng}
\end{figure}

Let $\mathcal{F}$ be the set of $52$ forbidden induced subgraphs of \hngsp found computationally on $6$ to $8$ vertices, shown in Figure \ref{fig:fis_for_1hng}. Let $\mathcal{F_S}$ denote those graphs in $\mathcal{F}$ containing no induced $C_5$, and $\mathcal{F_C}$ denote graphs in $\mathcal{F}$ containing an induced $C_5$. The subscript $S$ will denote that these are forbidden induced subgraphs of the sum-perfect graphs; see Theorem \ref{prop:sum_perfect}. The set $\mathcal{F}_S$ contains 24 graphs on 6 vertices and 2 graphs on 7 vertices. Twelve of the 24 graphs on 6 vertices are bipartite graphs with matching number $\nu = 3$. Note that a bipartite graph has $\nu = 3$ if and only if it has an induced perfect matching on $6$-vertices, if and only if it contains a $3K_2$ subgraph. The other 12 are complements of these graphs. The two graphs on 7 vertices are the sun graph with a pendant (degree one) vertex attached to a degree two vertex, and its complement.

The set $\mathcal{F_C}$ contains 22 graphs on 7 vertices and 4 graphs on 8 vertices, and is also by necessity closed under complementation.

From here, we prove that $\mathcal{F}$ is exactly the set of forbidden induced subgraphs for \hng. 

\begin{theorem}
\label{prop:FIS_for_1_hng}
    A graph $G$ is in \hngsp if and only if it contains no element of $\mathcal{F}$ as an induced subgraph. 
\end{theorem}

This work is made much simpler for graphs in \hngsp without $C_5$ by the fact that these graphs are exactly the sum-perfect graphs, and their forbidden induced subgraph characterization is known. We present the characterization by Litjens, Polak, and Sivaraman \cite{LiPoSi19}. 

\begin{theorem}
\label{prop:sum_perfect}
    A graph $G$ is sum-perfect if and only if it is $\mathcal{F_S} \cup \{C_5\}$-free. 
\end{theorem}

Our proof of Theorem \ref{prop:FIS_for_1_hng} will separately address graphs that do and do not contain an induced $C_5$. In anticipation of the former, we introduce notation and a proposition summarizing the contents of $\mathcal{F_C}$. If a graph $G$ contains no element of $\mathcal{F}$ as an induced subgraph, then the possible coexisting vertices outside $C_5$ and their adjacencies are limited.

Let $G$ be a graph containing vertices $C = \{c_1, c_2, c_3, c_4, c_5\} \in V(G)$ that induce a copy of $C_5$ along edges $v_1v_2, v_2v_3, v_3v_4, v_4v_5,$ and $v_1v_5$. (We will choose such a subgraph $C$ frequently throughout this paper, and always with the given edge set.) Let $D = V(G) - C$; we say that a vertex $v \in D$ is of type $N(v) \cap C$, according to its adjacencies in $C_5$, such as $\{c_1, c_2, c_4\}, \{c_3, c_5\}$, or even $\{\emptyset\}$. Multiple vertices in $D$ may be of the same type. The automorphism group of $C_5$ is the dihedral group on $5$ elements. Throughout, we will work up to symmetry of this automorphism group, and relabel $C$, up to symmetry, so that any chosen vertex or vertices are of lowest possible type. For instance, if we assume the existence of a vertex with one neighbor in $C$, we would call it of type $\{c_1\}$, up to symmetry. With a vertex of type $\{c_2, c_4\}$ or $\{c_2, c_5\}$, we would apply the requisite relabeling of $C$ to have a vertex of type $\{c_1, c_3\}$.

The following remark summarizes the pairs of vertices, up to symmetry, not in $C_5$ that can occur within $D$ without containing an induced graph from $\mathcal{F}$.   

\begin{prop}
\label{prop:allowable_vertex_list}
    Let $G$ be a graph containing no element of $\mathcal{F}$ as an induced subgraph, and let $C = \{c_1, c_2, c_3, c_4, c_5\} \in V(G)$ induce a copy of $C_5$ with edges $c_1c_2,$ $ c_2c_3, $ $c_3c_4, $ $c_4c_5,$ and $c_1c_5$. Let $D = V(G) - C$ and let $v \in D$.
    
    If $v$ is of type $\{\emptyset\}$, then for all $w \in D$, either $w$ is adjacent to $v$ and of types $\{c_1, c_2, c_3\},$$\{c_1, c_2, c_3, c_4\}$ or $\{c_1, c_2, c_3, c_4, c_5\}$, up to symmetry, or $w$ is non-adjacent to $v$ and of any type.
    
    If $v$ is of type $\{c_1\}$, up to symmetry, then for all $w \in D$, either $w$ is adjacent to $v$ and of type $\{c_3, c_4\}$ or $\{c_1, c_2, c_3, c_4, c_5\}$, or $w$ is non-adjacent to $v$ and of type $\{\emptyset\},\{c_1\},$ $ \{c_1, c_2\},$ $ \{c_1, c_3\}, $ $\{c_1, c_4\},$ $ \{c_1, c_5\},$ $ \{c_3, c_4\},$ $ \{c_1, c_2, c_4\}, $ $\{c_1, c_2, c_5\}, $ $\{c_1, c_3, c_4\}, $ $\{c_1, c_3, c_5\}$, or $\{c_1, c_2, c_3, c_4, c_5\}$. 
    
    If $v$ is of type $\{c_1, c_2\}$, up to symmetry, then for all $w \in D$, either $w$ is adjacent to $v$ and of type $\{c_4\},$ $\{c_1, c_2\}, $ $\{c_1, c_2, c_3\},$ $\{c_1, c_2, c_4\}, $ $\{c_1, c_2, c_5\},$ $\{c_1, c_2, c_3, c_4\},$ $\{c_1, c_2, c_3, c_5\},$ $\{c_1, c_2, c_4, c_5\}$ or $\{c_1, c_2, c_3, c_4, c_5\}$, or $w$ is non-adjacent to $v$ and of type $\{\emptyset\},$ $\{c_1\},$ $\{c_4\},$ $\{c_1, c_2\},$$\{c_1, c_4\},$ $\{c_2, c_4\}$, or $\{c_1, c_2, c_4\}$. 
    
    If $v$ is of type $\{c_1, c_3\}$, up to symmetry, then for all $w \in D$, either $w$ is adjacent to $v$ and of type $\{c_1, c_3, c_4, c_5\}$ or $\{c_1, c_2, c_3, c_4, c_5\}$, or $w$ is non-adjacent to $v$ and of type $\{\emptyset\},$ $\{c_1\},$ $\{c_3\},$ $\{c_1, c_3\},$ $\{c_1, c_4\},$ $\{c_1, c_5\},$ $\{c_3, c_4\},$ $\{c_3, c_5\},$ $\{c_1, c_3, c_4\},$ $\{c_1, c_3, c_5\}$, or $\{c_1, c_2, c_3, c_4, c_5\}$. 
    
    If $v$ is of type $\{c_1, c_2, c_3\}$, up to symmetry, then for all $w \in D$, either $w$ is adjacent to $v$ and of type $\{\emptyset\},$ $\{c_1, c_2\},$ $\{c_2, c_3\},$ $\{c_1, c_2, c_3\},$ $\{c_1, c_2, c_4\},$ $\{c_1, c_2, c_5\},$ $\{c_2, c_3, c_4\},$ $\{c_2, c_3, c_5\},$ $\{c_1, c_2, c_3, c_4\},$ $\{c_1, c_2, c_3, c_5\}$ or $\{c_1, c_2, c_3, c_4, c_5\}$, or $w$ is non-adjacent to $v$ and of type $\{\emptyset\}$ or $\{c_2\}$. 
    
    If $v$ is of type $\{c_1, c_2, c_4\}$, up to symmetry, then for all $w \in D$, either $w$ is adjacent to $v$ and of type $\{c_1, c_2\},$ $\{c_1, c_2, c_3\},$ $\{c_1, c_2, c_4\},$ $\{c_1, c_2, c_5\},$ $\{c_1, c_2, c_3, c_4\},$ $\{c_1, c_2, c_3, c_5\},$ $\{c_1, c_2, c_4, c_5\}$ or $\{c_1, c_2, c_3, c_4, c_5\}$, or $w$ is non-adjacent to $v$ and of type $\{\emptyset\},$ $\{c_1\},$ $\{c_2\},$ $\{c_4\},$ $\{c_1, c_2\},$ $\{c_1, c_4\},$ $\{c_2, c_4\},$ $\{c_1, c_2, c_4\}$, or $\{c_1, c_2, c_3, c_5\}$. 
    
    If $v$ is of type $\{c_1, c_2, c_3, c_4\}$, up to symmetry, then for all $w \in D$, either $w$ is adjacent to $v$ and of type $\{c_1, c_2\},$ $\{c_1, c_4\},$ $\{c_2, c_3\},$ $\{c_3, c_4\},$ $\{c_1, c_2, c_3\},$ $\{c_2, c_3, c_4\},$ $\{c_1, c_2, c_4\},$ $\{c_1, c_3, c_4\},$ $\{c_2, c_3, c_5\},$ $\{c_1, c_2, c_3, c_4\}$ or $\{c_1, c_2, c_3, c_4, c_5\}$, or $w$ is non-adjacent to $v$ and of type $\{\emptyset\}$ or $\{c_2, c_3, c_5\}$. 
    
    If $v$ is of type $\{c_1, c_2, c_3, c_4, c_5\}$, then for all $w \in D$, either $w$ is adjacent to $v$ and of any type, or $w$ is non-adjacent to $v$ and of type $\{\emptyset\},$ $\{c_1\},$ $\{c_2\},$ $\{c_3\},$ $\{c_4\},$ $\{c_5\},$ $\{c_1, c_3\},$ $\{c_1, c_4\},$ $\{c_2, c_4\},$ $\{c_2, c_5\}$, or $\{c_3, c_5\}$. 
\end{prop}
\begin{proof}
    It is straightforward to verify the result (by hand or computer). For $v$ of any type out of $\{\emptyset\},$ $\{c_1\},$ $\{c_1, c_2\},$ $\{c_1, c_3\},$ $\{c_1, c_2, c_3\},$ $\{c_1, c_2, c_4\},$ $\{c_1, c_2, c_3, c_4\},$ or $\{c_1, c_2, c_3, c_4, c_5\}$, if $w$ is of a listed type and adjacency to $v$, then the subgraph induced on vertices $\{v, w, c_1, c_2, c_3, c_4, c_5\}$ is neither isomorphic to nor contains an element of $\mathcal{F}$. If $w$ is not of a listed type and adjacency to $v$, then the subgraph induced on vertices $\{v, w, c_1, c_2, c_3, c_4, c_5\}$ is isomorphic to or contains an element of $\mathcal{F}$. 
\end{proof}

We will use this result profusely throughout our proof of Theorem \ref{prop:FIS_for_1_hng}.  

\begin{proof}
The reverse direction is straightforward to check: for each of the $52$ graphs in $\mathcal{F}$, each has $\chi(G) + \theta(G) = |G|-1$, but for all proper induced subgraphs $H$, $\chi(H) + \theta(H) \ge |H|$. 

For the forward direction, let $G$ be a graph containing none of the graphs in $\mathcal{F}$ as an induced subgraph. We consider separately graphs that do and do not contain an induced subgraph isomorphic to $C_5$. Suppose, first, that $G$ does not contain $C_5$ as an induced subgraph. Hence by Theorem \ref{prop:sum_perfect}, $G$ is sum-perfect, so $\omega(G) + \alpha(G) \ge n$. Since $\omega(G) \le \chi(G)$ and $\alpha(G) \le \theta(G)$, it follows that $chi(G) + \theta(G) \ge n$ and therefore $G \in$ \hng. 

Otherwise, let $G$ contain an induced copy of $C_5$. Suppose that vertices $\{c_1, c_2, c_3, c_4, c_5\} \in V(G)$ induce a copy of $C_5$ with edges $v_1v_2, v_2v_3, v_3v_4, v_4v_5,$ and $v_1v_5$. Let $C = \{c_1,c_2,c_3,c_4,c_5\}$ and $D = V(G) - C$. In any graph, isolated vertices are $\theta$-distinct and dominating vertices are $\chi$-distinct. Assume by Proposition \ref{prop:vertex_deletion_mfis} that $G$ has neither isolated nor dominating vertices. It is also straightforward to check that if $|D| \le 2$, then $G \in$ \hng, so we assume $|D| \ge 3$. 

We will exploit Proposition \ref{prop:allowable_vertex_list} to identify possible vertex types in $G$. Then, we show either (i), there exists a $\theta$-distinct vertex in $G$, or (ii), there exists some $v \in V(G)$ whose neighborhood is isolated and no vertex in it is $\theta$-distinct in $G - \{v\}$. If (i) holds, $G$ cannot be a forbidden induced subgraph by Proposition \ref{prop:vertex_deletion_mfis}, and if (ii) holds, Corollary \ref{prop:theta_distinct_mfis} implies the same. The proof consists of four cases: first, where there exists a vertex of type $\{\emptyset\}$. Since the forbidden induced subgraphs of \hngsp are closed under complementation (see Remark \ref{prop:FIS_closed_complement}), we can assume for subsequent cases that no vertex is of type $\{\emptyset\}$ or $\{c_1, c_2, c_3, c_4, c_5\}$. The second case is where there exists a vertex of type $\{c_1\}$, up to symmetry. Again, for subsequent cases, we assume there are no vertices of types $\{c_1\},$ $\{c_2\},$ $\{c_3\},$ $\{c_4\},$ $\{c_5\},$ $\{c_1, c_2, c_3, c_4\},$ $\{c_1, c_2, c_3, c_5\},$ $\{c_1, c_2, c_4, c_5\},$ $\{c_1, c_3, c_4, c_5\},$ or $\{c_2, c_3, c_4, c_5\}$. The third case is where there exists a vertex of type $\{c_1, c_3\}$, up to symmetry. For the fourth case, we assume there are only vertices of types $\{c_1, c_2\},$ $\{c_1, c_5\},$ $\{c_2, c_3\},$ $\{c_3, c_4\},$ $\{c_4, c_5\},$ $\{c_1, c_2, c_4\},$ $\{c_1, c_3, c_4\},$ $\{c_1, c_3, c_5\},$ $\{c_2, c_3, c_5\},$ or $\{c_2, c_4, c_5\}$, and suppose up to symmetry that there exists a vertex of type $\{c_1, c_2\}$. 

{\textbf {Case 1: $D$ contains a vertex of type $\{\emptyset\}$.} }

Let $E$ be the set of vertices of type $\{\emptyset\}$ and let $N(E)$ be the set of vertices with a neighbor in $E$. By Proposition \ref{prop:allowable_vertex_list}, $E$ induces a stable set, and by assumption, both $E$ and $N(E)$ are nonempty (else a vertex in $E$ would be isolated). By Proposition \ref{prop:allowable_vertex_list}, $N(E)$ can contain vertices of types $\{c_1, c_2, c_3\}, \{c_1, c_2, c_3, c_4\}$, and $\{c_1, c_2, c_3, c_4, c_5\}$, up to symmetry. Any vertices of these types must be adjacent by Proposition \ref{prop:allowable_vertex_list}, so $N(E)$ induces a clique. 

We now show $N(E)$ is complete to $B = D - E - N(E)$. First, let $w \in N(E)$ be a vertex of type $\{c_1, c_2, c_3\}$, up to symmetry. Proposition \ref{prop:allowable_vertex_list} implies $w$ is complete to $B$ except for vertices of type $\{c_2\}$. However, if there exists a vertex of type $\{c_2\}$, $G$ contains $F_2$ as an induced subgraph. By Proposition \ref{prop:allowable_vertex_list}, $D$ contains no vertices of types $\{c_1\},\{c_3\},\{c_4\},\{c_5\},\{c_1, c_3\},\{c_1, c_4\},\{c_2, c_4\},\{c_2, c_5\},$ or $\{c_3, c_5\}$, so any vertex of type $\{c_1, c_2, c_3, c_4, c_5\}$ in $N(E)$ is also complete to $B$. 

Second, let $w \in N(E)$ be of type $\{c_1, c_2, c_3, c_4\}$. Proposition \ref{prop:allowable_vertex_list} implies $w$ is complete to $B$ except for possibly vertices of type $\{c_2, c_3, c_5\}$. However, given a vertex $x$ of type $\{c_2, c_3, c_5\}$ not adjacent to $w$, $G$ contains $F_{24}$ as an induced subgraph. Furthermore, by Proposition \ref{prop:allowable_vertex_list}, $D$ contains no vertices of types $\{c_1\},$ $\{c_2\}$ $\{c_3\},$ $\{c_4\},$ $\{c_5\},$ $\{c_1, c_3\},$ $\{c_2, c_4\},$ $\{c_2, c_5\},$ or $\{c_3, c_5\}$. More subtly, $D$ also does not contain a vertex of type $\{c_1, c_4\}$ that is not adjacent to some vertex of type $\{c_1, c_2, c_3, c_4, c_5\}$, else $G$ contains an induced $\overline{F_2}$. Thus any vertex of type $\{c_1, c_2, c_3, c_4, c_5\}$ in $N(E)$ is complete to $B$. 

The third and final option is that $N(E)$ contains only vertices of type $\{c_1, c_2, c_3, c_4, c_5\}$. By Proposition \ref{prop:allowable_vertex_list}, any vertex of type $\{c_1, c_2, c_3, c_4, c_5\}$ is adjacent to all other vertices in $B$ except possibly those of type $\{c_1\}$ and $\{c_1, c_3\}$, up to symmetry. If there exists a vertex of type $\{c_1, c_2, c_3, c_4, c_5\}$ that is non-adjacent to a vertex of type $\{c_1\}$ in $D$, up to symmetry, then $\overline{G}$ contains a vertex of type $\{\emptyset\}$ with a neighbor of type $\{c_1, c_2, c_3\}$. If there exists a vertex of type $\{c_1, c_2, c_3, c_4, c_5\}$ that is non-adjacent to a vertex of type $\{c_1, c_3\}$ in $D$, up to symmetry, then $\overline{G}$ contains a vertex of type $\{\emptyset\}$ with a neighbor of type $\{c_1, c_2, c_3, c_4\}$. Since \hngsp is closed under complementation, we can assume $N(E)$ contains a vertex of type $\{c_1, c_2, c_3\}$ or $\{c_1, c_2, c_3, c_4\}$, both resolved above. Hence we can proceed with the assumption that $N(E)$ is complete to $B$. 

Fix a minimum proper clique covering of $C \cup B$. We can extend this to a minimum proper clique covering of $C \cup B \cup N(E)$ by adding $N(E)$ to the clique containing $c_2$. Since $N(E)$ is complete to $B$, the clique covering is proper. This extends to a proper clique covering of $G$ by giving each vertex in $E$ a distinct clique, so $\theta(G) \le \theta(C \cup B) + |E|$. Since $E$ induces a stable set and is isolated to $C \cup B$, $\theta(G) \ge \theta(C \cup B) + |E|$. We conclude the specified clique covering is minimum. Since it is $v$-distinct for any $v \in E$, Proposition \ref{prop:vertex_deletion_mfis} implies that $G$ is not a forbidden induced subgraph of \hng.

\textbf{Case $2$: $D$ contains a vertex of type $\{c_1\}$.}

Suppose that $D$ contains no vertices of type $\{\emptyset\}$ or $\{c_1, c_2, c_3, c_4, c_5\}$, but does contain a vertex of type $\{c_1\}$, up to symmetry. By Proposition \ref{prop:allowable_vertex_list}, $D$ may also contain vertices of type $\{c_1\}, $ $\{c_1, c_2\},$ $ \{c_1, c_3\},$ $ \{c_1, c_4\},$ $ \{c_1, c_5\},$ $ \{c_3, c_4\}, $ $\{c_1, c_2, c_4\},$ $ \{c_1, c_2, c_5\},$ $ \{c_1, c_3, c_4\},$ and $\{c_1, c_3, c_5\}$. Of these, vertices of type $\{c_1\}$ may only be adjacent to vertices of type $\{c_3, c_4\}$. 

Suppose there exists $v$ of type $\{c_1\}$ not adjacent to any vertex of type $\{c_3, c_4\}$, so $N(v) = \{c_1\}$ is a stable set. Suppose for a contradiction there exists a $c_1$-distinct clique covering of $H = G - \{v\}$. Since $N(c_2) \subseteq N(c_1) - \{c_3\}$ with the exception of $c_3$, it follows that $c_2$ and $c_3$ share a clique. By an analogous argument, $c_4$ and $c_5$ share a clique, and so no vertex in $H - C_5$ shares a clique with any vertex in $C_5$. By assumption, $|D| \ge 3$, so $|H - C_5|$ is nonempty. Choose a clique $K$ in $H - C_5$; there are three possibilities: 
\begin{enumerate}
    \item [(i)] $K$ has no vertex of type $\{c_3, c_4\}$. 
    If so, $c_1$ is complete to $K$ and can be added to it, contradicting the clique covering's minimality. 
    \item [(ii)] $K$ has a vertex of type $\{c_3, c_4\}$ but no vertex of type $\{c_1\}$. If so, $K$ contains only vertices of types $\{c_3, c_4\}$ and $\{c_1, c_3, c_4\}$, so we add $c_3$ to $K$ and combine $c_1$ and $c_2$ into one clique, contradicting the assumption that the clique covering is minimum.
    \item [(iii)] $K$ has a vertex of type $\{c_3, c_4\}$ and a vertex of type $\{c_1\}$. These two vertices, together with $C_5$ and $v$, induce $F_{25}$ or $F_{26}$ in $G$, for a contradiction. 
\end{enumerate}  
Instead, we conclude $H$ has no $c_1$-distinct clique covering, and Corollary \ref{prop:theta_distinct_mfis} implies $G$ is not a forbidden induced subgraph of \hng. 

Otherwise, every vertex of type $\{c_1\}$ is adjacent to a vertex of type $\{c_3, c_4\}$. If there exist two vertices of type $\{c_1\}$, they cannot be adjacent by Proposition \ref{prop:allowable_vertex_list}. If they are adjacent to a common vertex of type $\{c_3, c_4\}$, $G$ contains an induced $F_{26}$, and if they are adjacent to different vertices of type $\{c_3, c_4\}$, $G$ contains an induced $F_{25}$. Thus $G$ has exactly one vertex $v$ of type $\{c_1\}$. If there exist two vertices of type $\{c_3, c_4\}$ adjacent to $v$, then $G$ contains either $F_3$ or $F_4$ as an induced subgraph. Hence $v$ is adjacent to exactly one vertex $w$ of type $\{c_3, c_4\}$. The vertex $w$ is adjacent to all other vertices of type $\{c_3, c_4\}$, else $G$ contains an induced $F_3$, and to all vertices of type $\{c_1, c_3, c_4\}$, else $G$ contains an induced $F_{17}$. Hence $w$ is complete to $D$. 

Here, $N(v) = \{c_1, w\}$, which induces a stable set. We claim both $c_1$ and $w$ are not $\theta$-distinct in $H = G - \{v\}$. Suppose for a contradiction there exists a $w$-distinct clique covering of $H$. Thus $D - \{v\}$ shares a clique with $c_1$, $c_2$ and $c_3$ share a clique, and $c_4$ and $c_5$ share a clique. However, the clique covering is not minimum: put $c_1$ with $c_2$ and add $c_3$ and $w$ to $D - \{v\}$ to produce a strictly smaller clique covering of $H$. The covering is proper since $w$ is complete to $D$. Suppose for a contradiction there exists a $c_1$-distinct clique covering of $H$. Thus $c_2$ and $c_3$ share a clique in that covering, and $c_4$ and $c_5$ share a clique. Put $c_3$ in $w$'s clique and add $c_1$ to $c_2$'s clique to produce a proper, strictly smaller clique covering. Instead, we conclude $H$ has no $c_1$-distinct or $w$-distinct clique covering, and Corollary \ref{prop:theta_distinct_mfis} implies $G$ is not a forbidden induced subgraph of \hng. 

{\textbf {Case 3: $D$ contains a vertex of type $\{c_1, c_3\}$.} }

By ruling out the previous cases, $D$ contains no vertices with zero, one, four, or five neighbors in $C$, but does contain a vertex $v$ of type $\{c_1, c_3\}$, up to symmetry and complementation. By Proposition \ref{prop:allowable_vertex_list}, $D$ may also contain vertices of type $\{c_1, c_3\}, \{c_1, c_4\}, \{c_1, c_5\}, \{c_3, c_4\}, \{c_3, c_5\}, \{c_1, c_3, c_4\},$ and $\{c_1, c_3, c_5\}$. Any vertices of type $\{c_1, c_3\}$, $\{c_1, c_4\}$, or $\{c_3, c_5\}$ are isolated in $D$. However, $D$ contains vertices of type $\{c_1, c_4\}, \{c_3, c_4\}$, or $\{c_1, c_3, c_4\}$ if and only if it does not contain vertices of type $\{c_1, c_5\}, \{c_3, c_5\},$ or $\{c_1, c_3, c_5\}$. Up to symmetry, we assume the former. 

Here, $N(v) = \{c_1, c_3\}$, which induces a stable set. We claim both $c_1$ and $c_3$ are not $\theta$-distinct in $H = G - \{v\}$. Suppose for a contradiction there exists a $c_1$-distinct clique covering of $H$. Thus $c_2$ and $c_3$ share a clique, $c_4$ and $c_5$ share a clique, $D$ contains no vertex of type $\{c_1, c_4\}$. and no vertex in $H - C_5$ shares a clique with any vertex in $C_5$. Since $H - C_5$ is nonempty and complete to $c_3$, join $c_3$ to any clique in $H - C_5$ and combine $c_1$ and $c_2$ to produce a proper clique covering of $H$ with fewer cliques. If there exists a $c_3$-distinct clique covering of $H$, then $c_1$ and $c_2$ share a clique. Since $c_3$ cannot be added to any other clique, $H - C_5$ contains only vertices of type $\{c_1, c_4\}$ and by assumption must have at least $2$ such vertices. However, the proper clique covering of $H$ with $c_2$ and $c_3$ sharing a clique, $c_1$ and $c_5$ sharing a clique, and $c_4$ in a clique with some vertex of type $\{c_1, c_4\}$ contains fewer cliques. Instead, we conclude $H$ has no $c_1$-distinct or $c_3$-distinct clique covering, and Corollary \ref{prop:theta_distinct_mfis} implies $G$ is not a forbidden induced subgraph of \hng. 

{\textbf {Case 4: $D$ contains a vertex of type $\{c_1, c_2\}$.} }

Here, $D$ contains only vertices of types $\{c_1, c_2\},$ $ \{c_1, c_5\}, $ $\{c_2, c_3\}, $ $\{c_3, c_4\},$ $ \{c_4, c_5\},$ $ \{c_1, c_2, c_4\},$ $ \{c_1, c_3, c_4\}, $ $\{c_1, c_3, c_5\},$ $ \{c_2, c_3, c_5\}$ and $\{c_2, c_4, c_5\}$. Without loss of generality, suppose that $D$ contains a vertex $v$ of type $\{c_1, c_2\}$. Thus $D$ can also contain only vertices of types $\{c_1, c_2\}$ and $\{c_1, c_2, c_4\}$, which may or may not be adjacent to one another. With $D$ so limited, we move directly to a clique covering. 

The set $N(c_3) = \{c_2, c_4\}$ induces a stable set, and we claim that neither $c_2$ nor $c_4$ is $\theta$-distinct in $H = G - \{c_3\}$. Since $c_2$ is complete to $\{v\} \cup N(v)$, no clique covering of $H$ is $c_2$-distinct, else $c_2$ could be joined to $v$'s clique. If a clique covering is $c_4$-distinct, then $c_5$ and $c_1$ share a clique. Since $c_1$ is complete to $\{v\} \cup N(v)$, add $c_1$ to $v$'s clique and combine $c_4$ and $c_5$ into one clique to produce a proper clique covering of $H$ with fewer cliques, for a contradiction. We conclude by Corollary \ref{prop:theta_distinct_mfis} that $G$ is not a forbidden induced subgraph of \hng. 
\end{proof}

\section{Graphs in \hngsp are Apex-Perfect}

We prove that graphs in \hngsp are apex-perfect, that is, that the deletion of some vertex yields a perfect graph. This is sufficient to show that \hngsp is $\chi$-bounded by the function $f(x) = x+1$. 

\begin{theorem}
\label{prop:apex_perfect}
    If $G \in $ \hng, then there exists $v \in G$ such that $G - \{v\}$ is a perfect graph. 
\end{theorem}
\begin{proof}
    Let $G \in $ \hng. The graphs $C_6, \overline{C_6}, P_6,$ and $\overline{P_6}$ are elements of $\mathcal{F}$, so $G$ contains none of these, and hence $G$ contains no cycle $C_m$ or cycle-complement $\overline{C_m}$ for $m \ge 6$. Thus $G$ is perfect if and only if $G$ contains no induced $C_5$. From here, we consider the structure of induced copies of $C_5$ in $G$ to show that all induced copies of $C_5$ must intersect at a single vertex. 

    Assume that $G$ is imperfect (i.e., not a perfect graph) and let $C = \{c_1, c_2, c_3, c_4, c_5\} \in V(G)$ induce $C_5$ with edges $c_1c_2, c_2c_3, c_3c_4, c_4c_5,$ and $c_1c_5$. If $G$ has exactly one induced copy of $C_5$, then the theorem holds, so we assume another exists. In the remainder of the proof, we consider vertex types that can be included in another $C_5$, and show that in each case $G$ is apex-perfect. 

    The neighbor set of a vertex of type $\{\emptyset\}$ induces a clique, as shown in the proof of Case 1 of Theorem \ref{prop:FIS_for_1_hng}. Hence no vertices of type $\{\emptyset\}$ or, by complement, of type $\{c_1, c_2, c_3, c_4, c_5\}$ are elements of any induced $C_5$. 
    
    Suppose there exists a vertex $v$ of type $\{c_1\}$, up to symmetry and complementation, contained in an induced subgraph $X$ isomorphic to $C_5$. The only neighbors of $v$ in $G$ are $c_1$ and possibly some vertices of type $\{c_3, c_4\}$. Since our work in the proof of Case 2 of Theorem \ref{prop:FIS_for_1_hng} showed that a vertex of type $\{c_1\}$ cannot have two (nonadjacent) neighbors of type $\{c_3, c_4\}$, we conclude that $v$'s neighbors in $X$ are $c_1$ and a vertex $w$ of type $\{c_3, c_4\}$. If, for a contradiction, there exists a subset $Y \subseteq V(G)$ inducing $C_5$ with $c_1 \not \in Y$, then by Proposition \ref{prop:allowable_vertex_list}, the elements of $Y$ must be some of $c_2, c_3, c_4, c_5$ and vertices of types $\{c_3, c_4\}$ and $\{c_1, c_3, c_4\}$. Let $U$ denote the set of these possible elements. Since $c_2$ and $c_5$ have one neighbor in $U$, they cannot be in $Y$, and since $c_3$ and $c_4$ are adjacent to all but one element of $U$, they also cannot be in $Y$. Hence $Y$ is comprised entirely of vertices of types $\{c_3, c_4\}$ and $\{c_1, c_3, c_4\}$. Suppose up to complementation that at least three vertices of $Y$ are of type $\{c_3, c_4\}$. If three consecutive vertices of $Y$, say $v_1, v_2, v_3$ are of type $\{c_3, c_4\}$, then the induced subgraph of $G$ on vertex set $\{v_1, v_2, v_3, c_1, c_2, c_5\}$ is forbidden by Theorem \ref{prop:FIS_for_1_hng}. Otherwise, up to symmetry $v_1, v_2,$ and $v_4$ are of type $\{c_3, c_4\}$ and $v_3, v_5$ are of type $\{c_1, c_3, c_4\}$. Here, the induced subgraph on $\{v_2, v_3, v_4, c_1, c_2, c_5\}$ is forbidden by Theorem \ref{prop:FIS_for_1_hng}. In either case, we obtain the contradiction that $G \not \in$ \hng, so conclude instead that all induced copies of $C_5$ in $G$ contain vertex $c_1$. Thus its deletion renders the graph perfect.  Hence if $G$ contains a vertex of types $\{c_1\}$ or $\{c_1, c_2, c_3, c_4\}$, up to symmetry, contained in an induced $C_5$, the graph is apex-perfect. Assume for the remainder of the proof that $G$ contains no such vertex. 

    If a vertex $v$ of type $\{c_1, c_3\}$, up to symmetry and complementation, is in an induced $C_5$, then by Proposition \ref{prop:allowable_vertex_list}, its two neighbors in the $C_5$ must be $c_1$ and $c_3$. If, for a contradiction, there exists a subset $Y \subseteq V(G)$ inducing $C_5$ with $c_1 \not \in Y$, then by Proposition \ref{prop:allowable_vertex_list}, the elements of $Y$ must be some of $c_2, c_3, c_4, c_5$ and vertices of types $\{c_1, c_3\}$, $\{c_1, c_4\}$, $\{c_1, c_5\}$, $\{c_3, c_4\}$, $\{c_3, c_5\}$, $\{c_1, c_3, c_4\}$, and $\{c_1, c_3, c_5\}$. Let $U$ denote the set of these possible elements. Since $c_2$ has one neighbor in $U$, it cannot be in $Y$. Since we've shown that any vertex of type $\{c_x, c_{x+2}\}, $ with subscripts $\mod 5$ must have as its $C_5$ neighbors $c_x$ and $c_{x+2}$, it follows that $Y$ does not contain vertices of type $\{c_1, c_3\}$ or $\{c_1, c_4\}$. If $Y$ contains a vertex of type $\{c_3, c_5\}$, then $Y$ contains $c_3, c_5$, and vertices of types  Proposition \ref{prop:allowable_vertex_list} restricts $Y$ to containing $c_3, c_5$, and vertices of types $\{c_1, c_5\}$ and $\{c_1, c_3, c_5\}$. Since $c_5$ is adjacent to all but $c_3$, $Y$ does not induce $C_5$. Thus $Y$ does not contain a vertex of type $\{c_3, c_5\}$. If instead $Y$ contains a vertex of type $\{c_3, c_4\}$ or $\{c_1, c_3, c_4\}$, then $Y$ is restricted to containing $c_3, c_4$ and vertices of types $\{c_3, c_4\}$ or $\{c_1, c_3, c_4\}$, which was shown above to not be possible. Lastly, if $Y$ contains a vertex of type $\{c_1, c_5\}$ or $\{c_1, c_3, c_5\}$, then in order to induce $C_5$, $Y$ must contain two vertices $v_1, v_2$ of type $\{c_1, c_5\}$, two vertices $v_3, v_5$ of type $\{c_1, c_3, c_5\}$, and $c_3$. However, $G$ then contains the induced subgraph $\overline{F_{22}}$, for a contradiction. Instead, all induced copies of $C_5$ in $G$ contain vertex $c_1$, and its deletion renders the graph perfect.  Hence if $G$ contains a vertex of types $\{c_1, c_3\}$ or $\{c_1, c_2, c_3\}$, up to symmetry, contained in an induced $C_5$, the graph is apex-perfect. Assume for the remainder of the proof that $G$ contains no such vertex. 

    Lastly, if there is a vertex $v$ of type $\{c_1, c_2\}$, up to symmetry and complementation, in an induced $C_5$, then by Proposition \ref{prop:allowable_vertex_list}, the other vertices in that induced $C_5$ must be some of $c_1, c_2$, and vertices of types $\{c_1, c_2\}$ and $\{c_1, c_2, c_4\}$. Since $c_1$ and $c_2$ are adjacent to all but $c_4$, they do not occur in this $C_5$. As shown above, the $C_5$ cannot contain only vertices of types $\{c_1, c_2\}$ and $\{c_1, c_2, c_4\}$, so $c_4$ must be included. This $C_5$ then contains two vertices each of types $\{c_1, c_2\}$ and $\{c_1, c_2, c_4\}$, and $G$ contains $\overline{F_{22}}$ as an induced subgraph, for a contradiction.     
\end{proof}

It follows that graphs in \hngsp are $\chi$-bounded by the function $f(x) = x+1$, which is best possible. 

\begin{theorem}
\label{prop:chi_bound_Vizing}
    If $G \in$ \hng, then $\chi(G) \le \omega(G) + 1$. 
\end{theorem}
\begin{proof}
    Let $G \in $ \hng. By Theorem \ref{prop:apex_perfect}, there exists $v \in G$ such that $\chi(G - \{v\}) = \omega(G - \{v\})$. By Remark \ref{prop:vertex_deletion}, $\chi(G) \le \chi(G - \{v\}) + 1$. Either $\omega(G - \{v\}) = \omega(G)$ or $\omega(G - \{v\}) = \omega(G) - 1$. In the first case, $\chi(G) \le \chi(G - \{v\}) + 1 = \omega(G - \{v\}) + 1 = \omega(G) + 1$, and in the second case, $\chi(G) \le \chi(G - \{v\}) + 1 = \omega(G - \{v\}) + 1 = \omega(G)$. Both yield the desired result. 
\end{proof}

\section{The Intersection of \hngsp with Several Common Graph Classes}

In this section we provide equivalent structural and forbidden (induced) subgraph characterizations of \hngsp intersected with line graphs, claw-free graphs, and triangle-free graphs. 

\begin{figure}[ht]
\centering
 \includegraphics[height=7cm]{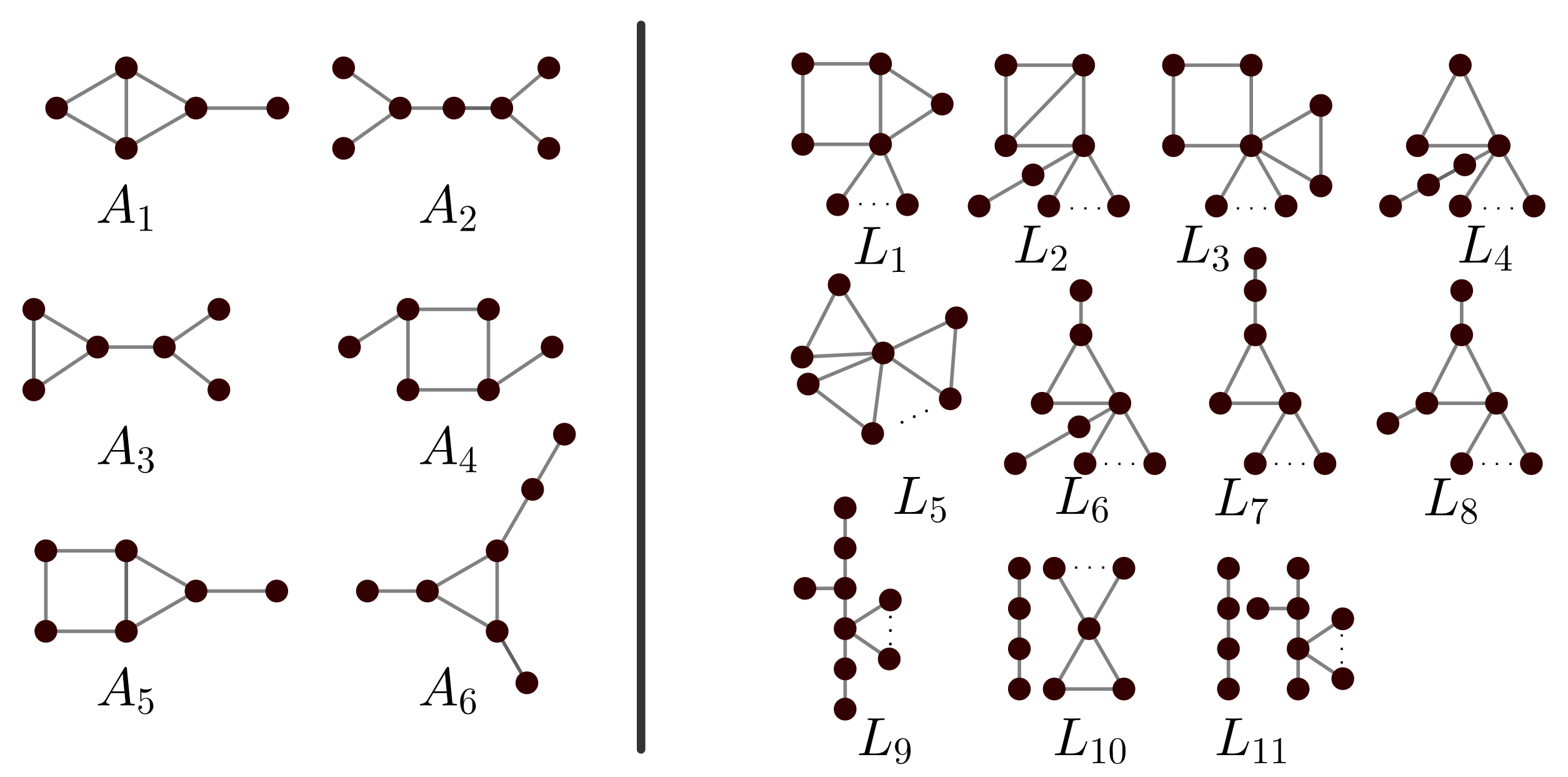}
 \caption{The set $\mathcal{A}$ comprises $A_1$ to $A_6$ as well as ten other small graphs. The set $\mathcal{L}$ comprises all graphs in families $L_1$ to $L_{11}$.}
  \label{fig:line_graphs}
\end{figure}

Beginning with line graphs, we define two sets that will provide our structural and forbidden subgraph characterizations. Given graphs $A_1$ to $A_6$ shown in Figure \ref{fig:line_graphs}, let $\mathcal{A} = \{A_1, \ldots, A_6, 3P_3, P_5 + P_3, P_7, C_6, K_4, C_4+P_3, 2K_3, 2K_{1,3}, K_{2,3}, K_3+K_{1,3}\}$. Given families $L_1$ to $L_{11}$ in Figure \ref{fig:line_graphs}, let $\mathcal{L}$ be the set of all graphs in one of families $L_1$ to $L_{11}$. We characterize the set of graphs whose line graphs are in \hngsp as the set of graphs containing none of $\mathcal{A}$ as a (not necessarily induced) subgraph, and equivalently as the set of graphs contained in an element of $\mathcal{L}$. 

\begin{theorem}
\label{prop:hng_line_graph}
Given a graph $G$, the following are equivalent: 
\begin{enumerate}
    \item[(i)] $L(G) \in $ \hng. 
    \item[(ii)] $G$ contains none of $\mathcal{A}$ as a subgraph. 
    \item[(iii)] $G$ is a subgraph of a graph in $\mathcal{L}$. 
\end{enumerate}
\end{theorem}
\begin{proof}
    (i) $\rightarrow$ (ii) 
    It is straightforward to verify that for all $A \in \mathcal{A}$, the line graph $L(A) \in \mathcal{F}$. Thus if a graph $G$ contains $A$ as a subgraph, $L(G)$ contains the forbidden $L(A)$ as an induced subgraph, and by Theorem \ref{prop:FIS_for_1_hng}, $L(G) \not \in$ \hng. 
    
    (ii) $\rightarrow$ (iii) 
    Let $G$ be a graph containing none of $\mathcal{A}$ as a subgraph. Since $C_6, P_7 \in \mathcal{A}$, it follows that $G$'s largest cycle has $5$ or fewer vertices, should one exist. It is clear that $G$ can contain arbitrarily many connected components isomorphic to $K_2$, and so we restrict the subsequent analysis to edgewise nontrivial components. We split into cases according to the number of components containing at least two edges and the size of the largest cycle. 

    Case 1: The graph $G$ has one edgewise nontrivial component, which contains a $5$-cycle. Any added vertex with two neighbors in the $5$-cycle either produces $C_6 \in \mathcal{A}$ as a subgraph if the neighbors are adjacent or $A_4 \in \mathcal{A}$ if the neighbors are not adjacent. Hence any vertex beyond the $5$-cycle must have at most one neighbor in the $5$-cycle. If two different vertices in the $5$-cycle have added neighbors, then $G$ contains either $P_7 \in \mathcal{A}$ (if the vertices are adjacent) or $A_2 \in \mathcal{A}$ (if not), so at most one vertex in the $5$-cycle has an added neighbor. If there exists an edge $e$ not incident to the $5$-cycle, then $e$ together with its path to the $5$-cycle and the $5$-cycle itself contain $P_7 \in \mathcal{A}$, so the only possible addition to the $5$-cycle is a set of pendant vertices attached to one vertex. If the $5$-cycle has two chords, then whether or not they are crossing chords, $G$ contains $A_1 \in \mathcal{A}$ as a subgraph. It follows that $G$ has at most one chord. A chord can be added if and only if it is incident to the vertex with pendant vertices, else $G$ contains $A_2$ or $A_4 \in \mathcal{A}$. It can be added freely and we conclude that $G$ is a subgraph of a graph in family $L_1$. 

    Case 2: The graph $G$ has one edgewise nontrivial component, which contains no $5$-cycle but does contain a $4$-cycle. Any added vertex with two neighbors in the $4$-cycle either produces $C_5$ (handled in Case 1) as a subgraph if the neighbors are adjacent or $K_{2,3} \in \mathcal{A}$ if the neighbors are not adjacent. Hence any added vertex has at most one neighbor in the $4$-cycle. If two different vertices in the $4$-cycle have added neighbors, then those two vertices must be adjacent, else $G$ contains $A_4 \in \mathcal{A}$. If two adjacent vertices in the $4$-cycle each have two added neighbors, then $G$ contains $2K_{1,3} \in \mathcal{A}$, so instead only one vertex in the $4$-cycle can have more than one additional neighbor. If there exists an edge $e$ not incident to the $4$-cycle and two vertices in the $4$-cycle have additional neighbors, then $e$ together with its path to the $4$-cycle, the $4$-cycle itself, and the additional neighbor contain $P_7 \in \mathcal{A}$. Hence either all additions to the $4$-cycle are pendant vertices and all but one are appended to the same vertex, creating a subgraph of a graph in family $L_1$, or all additions to the $4$-cycle are appended to the same vertex, which we call $v$. 

    However, we are limited in the subgraphs we can append to $v$: the set of vertices not included in the $4$-cycle cannot contain $P_3$, else $G$ contains $C_4+P_3 \in \mathcal{A}$. Hence we are limited to appending isolated copies of $P_3, K_3$, and pendant vertices to $v$. With two copies of $P_3$ appended to $v$, $G$ contains $P_5 + P_3 \in \mathcal{A}$, so we can append either one $K_3$ or one $P_3$ (and unlimited pendant vertices) to $v$. If the $4$-cycle contains two chords, then $G$ contains $K_4 \in \mathcal{A}$, so it has at most one chord. This chord must be incident to $v$, else $G$ contains $A_1$. If there is a chord, then $K_3$ cannot be appended to $v$, else $G$ contains $A_3 \in \mathcal{A}$. With a chord, $G$ is thus a subgraph of a graph in family $L_2$, and without a chord, $G$ is a subgraph of a graph in family $L_3$. 

    Case 3: The graph $G$ has one edgewise nontrivial component, which contains $K_3$ but no larger cycle. Let $C = \{v_1, v_2, v_3\} \subseteq V(G)$ induce $K_3$. There is exactly one path from any other vertex to $C$, else $G$ contains $2K_3 \in \mathcal{A}$ or some larger cycle, addressed above. Hence the subgraphs appended to each vertex in $C$ are isolated from one another. If $P_4$ is appended with a vertex in $C$ as an endpoint and $G$ contains any other edge not incident to that vertex, then $G$ contains either $P_7$ or $P_5+P_3$. Hence if $P_4$ is appended to a vertex in $C$, $G$ is a subgraph of a graph in family $L_4$. 

    A claw cannot be appended at a spoke to a vertex in $C$, since this would create $A_3 \in \mathcal{A}$ as a subgraph. Hence the subgraph appended to any vertex of $C$ contains only isolated copies of $P_3, K_3$, and pendant vertices. If $K_3$ is appended to a vertex in $C$, then appending neighbors to any other vertex in (either copy of) $K_3$ yields $A_3$ as a subgraph. Any number of copies of $K_3$ can be appended at the same vertex, and we conclude that $G$ is a subgraph of a graph in family $L_5$. If two copies of $P_3$ are appended to $C$ at the same vertex, then no other vertex in $C$ has an additional neighbor, else $G$ contains $P_5 + P_3$. If two copies of $P_3$ are appended to different vertices, then $G$ contains $P_7$. If two vertices in $C$ each have two additional neighbors, then $G$ contains $A_2$. If one vertex has a copy of $P_3$ appended and the other two vertices in $C$ have additional neighbors, then $G$ contains $A_6 \in \mathcal{A}$. Thus with a $P_3$ appended to a vertex in $C$, either that vertex can have a cluster of pendant vertices and $G$ is a subgraph of a graph in $L_6$ or a different vertex has a cluster of pendant vertices and $G$ is a subgraph of a graph in $L_7$. With no $P_3$ appended to a vertex in $C$, $G$ is a subgraph of a graph in $L_8$. 

    Case 4: The graph $G$ has one edgewise nontrivial component, which contains no cycles. If $G$ contains $P_6$ on vertices $v_1, \ldots, v_6$ in that order, then it can have a cluster of pendant vertices appended to any non-endpoint vertex. No $P_3$ can be appended, else $G$ contains $P_7$ or $P_5+P_3$. Furthermore, two non-adjacent vertices cannot have appended vertices, else $G$ contains $A_2$. If $v_2$ has a pendant vertex, then $v_3$ cannot, else $G$ contains $P_5 + P_3$. If $v_3$ and $v_4$ each have two or more pendant vertices, then $G$ contains $2K_{1,3}$. Thus $G$ is isomorphic to a subgraph of $P_6$ together with a collection of pendant vertices to $v_2$, which is a subgraph of a graph in $L_7$; or $G$ has a collection of pendant vertices to $v_3$ and a pendant vertex attached to $v_4$ (up to symmetry), and hence is a subgraph of a graph in $L_9$. 

    If $G$ contains no induced $P_6$ but does contain an induced $P_5$ on vertices $v_1, \ldots, v_5$ in that order, then it can have adjoined copies of $P_3$, but only to $v_3$. The resulting graph, with a cluster of copies of $P_3$ adjoined to $v_3$, is a subgraph of a graph in $L_5$. Without an adjoined $P_3$, again $G$ can only have one cluster of pendant vertices, and only neighboring vertices of the $P_5$ can have pendant vertices. The result is a subgraph a graph in $L_9$. 

    Case 5: The graph $G$ has at least two edgewise nontrivial components. Since $G$ cannot contain $3P_3 \in \mathcal{A}$ as an induced subgraph, $G$ has at most two components with at least two edges. Components containing a single edge or vertex are trivial, and ignored henceforth. The nontrivial components cannot contain $C_4$ or $P_5$, since $C_4 + P_3$ and $P_5 + P_3$ are in $\mathcal{A}$. If one component has a triangle, there cannot be a copy of $P_3$ appended to one triangle vertex, additional neighbors appended to two triangle vertices, or a vertex of degree two or more appended to the triangle, since each would create a $C_4$ or $P_5$ subgraph. Hence a component with a triangle consists of, at most, a triangle with a cluster of pendant vertices adjoined to one vertex. The other component cannot contain a triangle or a claw since $2K_3, K_3 + K_{1,3} \in \mathcal{A}$, so $H$ must be a subgraph of $P_4$. Thus $G$ is a subgraph of a graph in $L_{10}$. Otherwise, both components are trees. If both components are subgraphs of $P_4$, $G$ is a subgraph of a graph in $L_{10}$. Otherwise, exactly one component contains a claw. This component contains $P_4$ with a pendant vertex attached to one of the center vertices. Since both center vertices cannot have two or more pendant vertices simultaneously without producing a $2K_{1,3}$ subgraph, the component is a subgraph of $P_4$ together with a cluster of pendant vertices at one center vertex and a single pendant vertex attached to the other center vertex. The other nontrivial component is again a subgraph of $P_4$, and we conclude that $G$ is a subgraph of a graph in $L_{11}$. 

    (iii) $\rightarrow$ (i) 
    Compute the line graph of each graph in $\mathcal{L}$. By Theorem \ref{prop:FIS_for_1_hng}, each contains no forbidden induced subgraph, so is in \hng. The same is true of its subgraphs. 
\end{proof}

\begin{figure}[ht]
\centering
 \includegraphics[height=3cm]{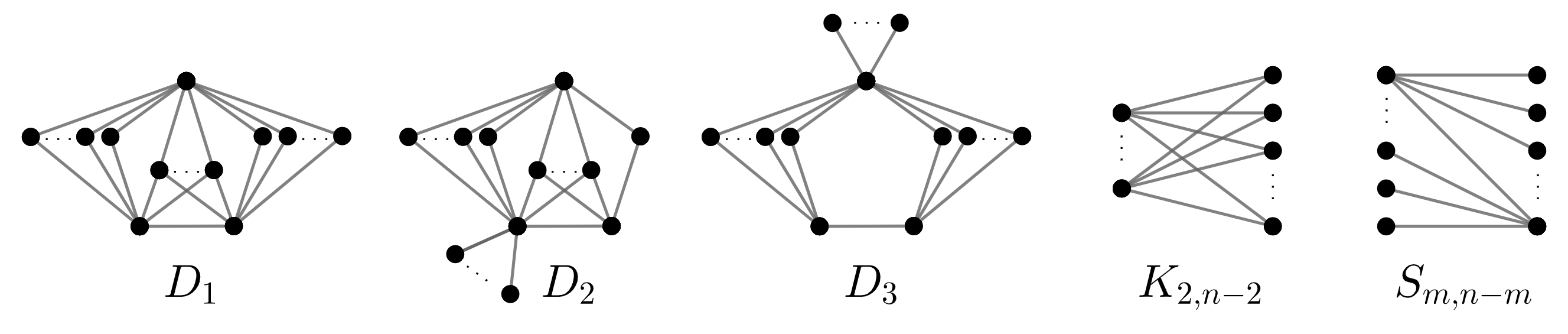}
 \caption{A claw-free graph in \hngsp containing an induced $C_5$ is an induced subgraph of a graph in family $\overline{D_1}, \overline{D_2}$, or $\overline{D_3}$, possibly with dominating vertices. A triangle-free graph in \hngsp is a subgraph of $K_{2, n-2}$ or $S_{m, n-m}$ for some $3 \le m \le n$, or an induced subgraph of $D_3$, possibly with isolated vertices.}
  \label{fig:claw_triangle_free}
\end{figure}

In preparation for characterizing claw-free and triangle-free subsets of \hng, we show that three families of graphs built from a copy of $C_5$ are in \hng. The three families, $D_1, D_2, $ and $D_3$, are shown in Figure \ref{fig:claw_triangle_free}.  Each consists of an induced $C_5$ on vertices labeled $\{c_1, c_2, c_3, c_4, c_5\}$ together with a stable set of added vertices. In $D_1$, the added vertices may be of types $\{c_1, c_3\}, $ $ \{c_1, c_4\}, $ and  $\{c_1, c_3, c_4\}$. In $D_2$, the added vertices may be of types $\{c_4\}, $ $ \{c_1, c_4\}, $ and  $\{c_1, c_3, c_4\}$. In $D_3$, the added vertices may be of types $\{c_1\}, $ $ \{c_1, c_3\}, $ and  $\{c_1, c_4\}$. 

\begin{lemma}
\label{prop:three_hng_families}
    All graphs in families $D_1, D_2$, or $D_3$ are in \hng. 
\end{lemma}
\begin{proof}
    Let $G$ be a graph in one of these families, with a fixed induced $C_5$ on vertex set $C = \{c_1, c_2, c_3, c_4, c_5\}$ in that order. By Theorem \ref{fig:fis_for_1hng} the addition to $C$ of any three vertices from $D$ does not induce an element of $F_C$. The choice of $C_5$ is immaterial: all induced copies of $C_5$ contain $c_1, c_3, c_4$, either $c_2$ or a vertex of type $\{c_1, c_3\}$, and either $c_5$ or a vertex of type $\{c_1, c_4\}$. With respect to any of these induced $C_5$, the types of vertices in $D = V(G) - C$ remain constant. Therefore $G$ contains no element of $F_C$ on any vertex subset. 

    Graphs in $D_1$ and $D_2$ may contain induced triangles, but all triangles share two common vertices, $c_3$ and $c_4$. Graphs in $D_3$ contain no induced triangles. In either case, $G$ cannot contain an induced $\overline{F_1}, \ldots, \overline{F_{12}}, F_13, $ or $\overline{F_{13}}$. See Figure \ref{fig:fis_for_1hng} for these graphs.  

    If $G$ is an element of $D_1$ or $D_3$, then the largest stable set containing $c_1$ has size $2$. Hence if $G$ contains an induced $F_1, \ldots, F_{11}, $ or $F_{12}$, it cannot include $c_1$. Without $c_1$, however, every edge in $G$ is incident to $c_3$ or $c_4$, so $\nu(G - \{c_1\}) = 2$. Hence $G$ contains no induced $F_1, \ldots, F_{12}$. If $G$ is an element of $D_2$, an analogous argument holds with the exchange of $c_1$ and $c_4$. 

    Therefore $G$ contains no induced element of $\mathcal{F}$, and Theorem \ref{fig:fis_for_1hng} implies that $G \in $\hng. 
\end{proof}

Given graphs $F_i$, $1 \le i \le 26$, shown in Figure \ref{fig:fis_for_1hng}, let $\mathcal{B} = \{K_{1,3}, F_1, F_2, F_3, F_5, F_8, F_{13},\overline{F_{21}}\} \cup \{\overline{F_1}, \overline{F_2}, \ldots, \overline{F_{13}}\}$, as defined in Figure \ref{fig:fis_for_1hng}. We use $\mathcal{B}$ as well as Lemma \ref{prop:three_hng_families} to characterize the claw-free graphs in \hng. 

\begin{theorem}
\label{prop:hng_claw_free}
Given a graph $G$, the following are equivalent: 
\begin{enumerate}
    \item[(i)] $G \in $ \hngsp and $G$ is claw-free. 
    \item[(ii)] $G$ contains no element of $\mathcal{B}$ as an induced subgraph. 
    \item[(iii)] Either $G$ is perfect, claw-free, and in \hng; or $\overline{G}$ is a graph in family $D_1, D_2,$ or $D_3$, possibly with the addition of some dominating vertices.  
\end{enumerate}
\end{theorem}
\begin{proof}
    (i) $\rightarrow$ (iii) 
    Let $G \in$ \hngsp be claw-free. Assume that $G$ contains no isolated vertices. By Theorem \ref{prop:apex_perfect}, $G$ is perfect if and only if it is $C_5$-free. If so, we are done, and if not, suppose $G$ contains an induced $C_5$ on vertex set $C = \{c_1, c_2, c_3, c_4, c_5\}$ in that order. Let $D = V(G) - C$. Since $G$ is claw-free, $D$ cannot contain any vertex $v$ of type $\{c_1\}, \{c_1, c_3\}, $ or $\{c_1, c_3, c_4\}$, up to symmetry, else $\{c_1, c_2, c_5, v\}$ induces a claw. Additionally, $D$ contains no vertex $v$ of type $\{\emptyset\}$, since then by Proposition \ref{prop:allowable_vertex_list}, its neighbor set contains a vertex $w$ with two non-adjacent neighbors in $C$, inducing a claw. 

    By virtue of the remaining permissible types, Proposition \ref{prop:allowable_vertex_list} implies $D$ must induce a clique. If $D$ contains a vertex $v$ of type $\{c_1, c_2\}$, up to symmetry, then Proposition \ref{prop:allowable_vertex_list} allows $D$ contain vertices of types $\{c_1, c_2\}, $ $\{c_1, c_2, c_3\},$ $\{c_1, c_2, c_5\},$ $\{c_1, c_2, c_3, c_4\},$ $\{c_1, c_2, c_3, c_5\},$ $\{c_1, c_2, c_4, c_5\}$ or $\{c_1, c_2, c_3, c_4, c_5\}$. However, if $D$ contains a vertex $w$ of type $\{c_1, c_2, c_3, c_5\}$ or $\{c_1, c_2, c_3, c_4, c_5\},$ then $\{v, w, c_3, c_5\}$ induces a claw. Since $D$ cannot simultaneously contain vertices of types $\{c_1, c_2, c_5\} $ and $ \{c_1, c_2, c_3, c_4\}$ or $\{c_1, c_2, c_4, c_5\}$ and $ \{c_1, c_2, c_3, c_4\}$ by Proposition \ref{prop:allowable_vertex_list}, it follows that $D$ contains vertices of types from set $\{\{c_1, c_2\},$ $ \{c_1, c_2, c_3\}, $ $\{c_1, c_2, c_5\}\}$ or $\{\{c_1, c_2\},$ $ \{c_1, c_2, c_3\},$ $ \{c_1, c_2, c_3, c_4\}\}$. If the former, then $\overline{G}$ is in $D_1$, and if the latter, $\overline{G}$ is in $D_2$. If $D$ contains no vertices of type $\{c_1, c_2\}$ up to symmetry, then the vertices of $D$ are of types from the set $\{\{c_1, c_2, c_3\}, \{c_1, c_2, c_5\}, \{c_1, c_2, c_3, c_4\}, \{c_1, c_2, c_3, c_5\}, \{c_1, c_2, c_3, c_4, c_5\}\}$, up to symmetry. However, if $D$ contains a vertex of type $\{c_1, c_2, c_3, c_4\},$ it cannot contain vertices of types $\{c_1, c_2, c_5\}$ or $ \{c_1, c_2, c_3, c_5\}$. Up to symmetry, we conclude the vertices of $D$ are from set $\{\{c_1, c_2, c_3\}$, $ \{c_1, c_2, c_5\}$, $ \{c_1, c_2, c_3, c_5\}$, $ \{c_1, c_2, c_3, c_4, c_5\}\}$, and so $G$ is in family $D_3$. 

    (iii) $\rightarrow$ (ii) 
    Let $G$ be a graph such that $\overline{G}$ is in family $D_1, D_2$, or $D_3$, possibly with dominating vertices. Since dominating and isolated vertices do not affect inclusion in \hng, suppose $G$ has none. Lemma \ref{prop:three_hng_families} implies $G$ contains no element of $\mathcal{B}$ as an induced subgraph, except possibly the claw.  
    
    We now show that if $G$ contains an induced $C_5$, it is claw-free. Suppose $\overline{G}$ is in $D_1$ or $D_2$. Since $D$ induces a clique, the only stable set of at least three vertices is $\{v, c_3, c_5\}$, and no vertex is adjacent to all three of these vertices. Thus $G$ contains no induced claw. Otherwise, if $\overline{G}$ is in $D_3$, then since $D$ induces a clique, each vertex in $D$ has at most two non-neighbors, which are adjacent. Thus $\alpha(G) = 2$ and $G$ contains no induced claw. 
    
    (ii) $\rightarrow$ (i) 
    Since $\mathcal{B}$ is precisely the subset of graphs in $\mathcal{F}$ not containing an induced claw, together with the claw, the result follows by Theorem \ref{prop:FIS_for_1_hng}. 
\end{proof}

In anticipation of characterizing triangle-free graphs of \hng, we provide a lemma about the structure of bipartite graphs with no $3K_2$ subgraph. By Theorem \ref{prop:FIS_for_1_hng}, this result also characterizes bipartite graphs in \hng. Let $S_{m, n-m}$ denote the double star graph with center vertices of degree $m$ and $n-m$; see Figure \ref{fig:claw_triangle_free}. 

\begin{lemma}
\label{prop:bipartite_doublestar}
    A graph $G$ is bipartite and has no $3K_2$ subgraph if and only if $G$ is a subgraph of $K_{2,n-2}$ or $S_{m, n-m}$. 
\end{lemma}
\begin{proof}
    Suppose $G$ is a subgraph of a double star graph with non-pendant vertices $v$ and $w$; hence $G$ is bipartite. Since all edges of $G$ are incident to $v$ or $w$, $G$ does not contain $3K_2$ as a subgraph. 

    Alternatively, suppose $G$ is bipartite and does not contain $3K_2$ as a subgraph. If one side of any possible bipartition has two or fewer vertices, then $G$ is a subgraph of $K_{2,n-2}$. Otherwise, in any bipartition $V(G) = A \cup B$ of $G$, both sides have three or more vertices, and suppose for a contradiction that $G$ is not a subgraph of a double star graph. Thus two vertices in $A$ must have different neighbors in $B$; suppose $E(G)$ contains $a_1b_1$ and $a_2b_2$. If there exists an edge with endpoints outside $a_1, a_2, b_1, b_2$, $G$ must contain $3K_2$; else, all edges are incident to one of these vertices. If there exists $a_3$ adjacent to $b_1$ and $b_3$ adjacent to $a_1$, then $a_1b_3, a_2b_2, a_3b_1$ form $3K_2$, and analogously if there exists $a_4$ adjacent to $b_2$ and $b_4$ adjacent to $a_2$. Otherwise, without loss of generality, all vertices in $A$ besides $a_1$ and $a_2$ are adjacent at most to $b_1$, and all vertices in $B$ besides $b_1$ and $b_2$ are adjacent at most to $a_2$. This is a subgraph of a double star graph unless $E(G)$ contains $a_1b_2$. In that case, since $A$ and $B$ each contain a third non-isolated vertex $a_3$ and $b_3$, the edges $a_1b_2, a_2b_3, a_3b_1$ form $3K_2$.
\end{proof}

To end this section, we characterize the triangle-free graphs in \hng. The structural characterization in (iii) is illustrated in Figure \ref{fig:claw_triangle_free}. 

\begin{theorem}
\label{prop:hng_triangle_free}
Given a graph $G$, the following are equivalent: 
\begin{enumerate}
    \item[(i)] $G \in $ \hngsp and $\omega(G) = 2$. 
    \item[(ii)] $G$ contains no induced $K_3, F_1, \ldots, F_{12}$. 
    \item[(iii)] $G$ is either (1) a subgraph of $K_{2,n-2}$, (2) a subgraph of $S_{m, n-m}$, or (3) a graph in family $D_3$, possibly with isolated vertices. 
\end{enumerate}
\end{theorem}
\begin{proof}
    (i) $\rightarrow$ (iii) 
    Let $G \in$ \hngsp with $\omega(G) = 2$. By Theorem \ref{prop:apex_perfect}, $G$ is perfect if and only if it is $C_5$-free. If so, $\omega = \chi = 2$ and so $G$ is bipartite. By Theorem \ref{prop:FIS_for_1_hng}, $G$ has no $3K_2$ subgraph, and so by Lemma \ref{prop:bipartite_doublestar}, $G$ is a subgraph of $K_{2,n-2}$ or a double star graph. Otherwise, $G$ is imperfect and contains an induced $C_5$ on some vertex subset $\{c_1, c_2, c_3, c_4, c_5\}$. Any other vertices must be of types $\{\emptyset\}, \{c_1\}$, and $\{c_1, c_3\}$, up to symmetry, since other types would induce $K_3$. The result follows from Proposition \ref{prop:allowable_vertex_list}. 

    (iii) $\rightarrow$ (ii) 
    If $G$ is a subgraph of $K_{2,n-2}$ or the double star graph, it clearly is bipartite and contains no induced $K_3$, and Lemma \ref{prop:bipartite_doublestar} implies $\nu(G) \le 2$, so $G$ contains no induced $F_1, \ldots, F_{12}$. If $G$ is a graph in family $D_3$, possibly with isolated vertices, then $G$ contains no induced $K_3$. By Lemma \ref{prop:three_hng_families}, $G \in $ \hng. 
    
    (ii) $\rightarrow$ (i) 
    Clearly $\omega(G) = 2$. All forbidden induced subgraphs of \hng, as laid out in Theorem \ref{prop:FIS_for_1_hng}, are given in (ii) or contain $K_3$, so we conclude that $G \in$ \hng.
\end{proof}

\section{Optimization}

In this section we show that membership in \hngsp, and many of its graphs' key invariants, can be determined in polynomial time. These results largely follow from the forbidden induced subgraph characterization and work in Section 4 on apex-perfection. 

\begin{theorem}
\label{prop:HNG_polynomialtime}
    Graphs in \hngsp can be recognized in polynomial time, specifically $O(n^8)$. 
\end{theorem}
\begin{proof}
    It suffices to test whether the given graph contains one of the $52$ graphs in the list, each of which has at most $8$ vertices.
\end{proof}

Furthermore, the clique number and independence number of graphs in \hngsp can be determined in polynomial time. Where the graphs are perfect, the results are immediate from the result by Chudnovsky et al. \cite{Ch05} that perfect graphs can be identified in polynomial time and from Gr\"{o}tschel, Lov\'{a}sz, and Schrijver's earlier work \cite{GrLoSc84} that perfect graphs can be given a minimum proper coloring in polynomial time. For imperfect graphs, we begin with a lemma refining the proof of Theorem \ref{prop:apex_perfect} to add an additional requirement to the deleted vertex, then continue with the full result in Theorem \ref{prop:cliquenumber_polynomialtime}.

\begin{lemma}
\label{prop:apexperfect_cliquenumber}
    If $G \in $ \hngsp and is imperfect, then there exists a vertex $v$ such that $G - \{v\}$ is perfect and $\omega(G - \{v\}) = \omega(G)$. 
\end{lemma}
\begin{proof}
    Let $G \in$ \hngsp and suppose $G$ is imperfect. By Theorem \ref{prop:apex_perfect}, $G$ is apex-perfect, and contains at least one induced $C_5$ (note again that $G$ cannot contain any of the other forbidden induced subgraphs of perfect graphs). Let $C = \{c_1, c_2, c_3, c_4, c_5\}$ induce a copy of $C_5$ in $G$ with edge set $\{c_1c_2, c_2c_3, c_3c_4, c_4c_5, c_1c_5\}$ and let $D = V(G) - C$. If $D = \emptyset$, then $G$ is isomorphic to $C_5$ and $\omega(G) = \omega(G - \{v\})$ for all $v \in V(G)$. If $\omega(G) = 2$, then since $G$ contains $C_5$, $\omega(G - \{v\}) = 2$ for all $v \in V(G)$. Otherwise, suppose that $D$ is nonempty and $\omega(G) \ge 3$. 

    If $D$ only has vertices of types $\{c_1, c_2\}$ and $\{c_1, c_2, c_4\}$, up to symmetry, then by the proof of Theorem \ref{prop:apex_perfect}, $G$ cannot contain an induced $C_5$ other than $C$. Since $c_3$ is contained in no triangles, its deletion verifies the theorem. Assume thus going forward that $G$ contains an induced $C_5$ other than $C$. 

    Suppose next that $D$ may also have vertices of types $\{c_1, c_3\}$ and $\{c_1, c_2, c_3\}$, up to symmetry. If an induced $C_5$ contains a vertex of type $\{c_1, c_3\}$, then its neighbors in the copy of $C_5$ must be $c_1$ and $c_3$, by Proposition \ref{prop:allowable_vertex_list}. The deletion of either verifies apex-perfection, so delete whichever of the two is not contained in a unique maximum clique. If an induced $C_5$ contains a vertex of type $\{c_1, c_2, c_3\}$, then its non-neighbors must analogously be $c_4$ and $c_5$. By Proposition \ref{prop:allowable_vertex_list}, $D$ cannot contain a vertex adjacent to both, so each is in a clique of size at most two, which is non-maximum by assumption. Delete either to reach the desired conclusion. 

    Suppose next that $D$ may also have vertices of types $\{c_1\}$ and $\{c_1, c_2, c_3, c_4\}$, up to symmetry. If an induced $C_5$ contains neither, then the above argument holds to verify the theorem. If an induced $C_5$ contains a vertex of type $\{c_1\}$, then its neighbors in the $C_5$ must be $c_1$ and a vertex of type $\{c_3, c_4\}$. By Proposition \ref{prop:allowable_vertex_list}, $c_1$ may neighbor vertices of types $\{c_1\}, \{c_1, c_3\}, \{c_1, c_4\}$ and $\{c_1, c_3, c_4\}$. All of these are isolated to one another, $c_2$, and $c_5$, with the exception of vertices of type $\{c_1, c_3, c_4\}$, which may be adjacent to one another. Hence if $c_1$ is in a maximum clique, the clique must comprise $c_1$ and at least two adjacent vertices of type $\{c_1, c_3, c_4\}$. Replace $c_1$ with $c_3$ and $c_4$ to obtain a strictly larger clique, for a contradiction. We conclude that $c_1$ is in no maximum clique, so its deletion verifies apex-perfection, by the proof of Theorem \ref{prop:apex_perfect}, and $\omega(G - \{c_1\}) = \omega(G)$. 

    If an induced $C_5$ contains a vertex of type $\{c_1, c_2, c_3, c_4\}$, then its non-neighbors in the $C_5$ must be $c_5$ and a vertex of type $\{c_2, c_3, c_5\}$. Since this case is the complement of the case in the proof of Theorem \ref{prop:apex_perfect} where an induced $C_5$ contains a vertex of type $\{c_1\}$, it follows that $G - \{c_5\}$ is perfect. By Proposition \ref{prop:allowable_vertex_list}, $c_5$ is in a clique of size at least $3$ only if all other clique vertices are of type $\{c_2, c_3, c_5\}$, but these vertices together with $c_2$ and $c_3$ form a strictly larger clique. Hence $\omega(G - \{c_5\}) = \omega(G)$. 

    Lastly, suppose $D$ may also contain vertices of types $\{\emptyset\}$ and $\{c_1, c_2, c_3, c_4, c_5\}$. Neither can be in any induced copies of $C_5$. Any vertex of type $\{c_1, c_2, c_3, c_4, c_5\}$ is included in all maximum cliques in $G$ and its induced subgraphs $G - \{v\}$ for all other $v \in V(G)$, since its only possible non-neighbors, vertices of types $\{c_1\}$ or $\{c_2, c_5\}$, up to symmetry, are never included in any maximum cliques. (Any two adjacent neighbors of a vertex of these types are adjacent to both $c_3$ and $c_4$, which would produce a strictly larger clique). Hence a vertex of type $\{c_1, c_2, c_3, c_4, c_5\}$ cannot impact the relative size of a maximum clique, and we proceed assuming none exist in $G$. 

    Any vertex of type $\{\emptyset\}$ is included in no maximum cliques in $G$ or its induced subgraphs $G - \{v\}$ for all other $v \in V(G)$, since its neighbor set (if nonempty) forms a clique dominating two adjacent vertices of $C$. Hence a vertex of type $\{\emptyset\}$ cannot change the clique number in $G$, and we can assume no such vertex exists in $G$. 
\end{proof}

Note also the complementary result pertaining to $\alpha(G)$ holds since \hng, perfect graphs, and all processes used in the proof are invariant under complementation. 

\begin{theorem}
\label{prop:cliquenumber_polynomialtime}
    If $G \in $ \hng, then $\omega(G)$ and $\alpha(G)$ can be determined in polynomial time. 
\end{theorem}
\begin{proof}
    If $G$ is perfect, the result holds immediately, as shown in \cite{GrLoSc84}. Otherwise, assume $G$ is imperfect. We show $\omega(G)$ can be determined in polynomial time; by complementation, the same holds for $\alpha(G)$. As established in the proof of Theorem \ref{prop:apex_perfect}, $G$ contains an induced copy of $C_5$.     
    In $O(n^5)$ time, identify a subset of vertices $\{v_1, v_2, v_3, v_4, v_5\}$ inducing $C_5$. For each $v_i$, $1 \le i \le 5$, determine for which $G - \{v_i\}$ is perfect, which can be done in polynomial time by \cite{Ch05}, and if so compute $\omega_i = \omega(G - \{v_i\})$, which is again possible in polynomial time. In all cases, $\omega_i = \omega(G)$ or $\omega(G) - 1$. By Lemma \ref{prop:apexperfect_cliquenumber}, this must hold with equality in the case of at least one $\omega_i$, so we conclude that $\omega(G) = \max\{\omega_i\}$. 
\end{proof}

We now perform a similar analysis to exactly determine the chromatic number of a graph in \hng, and a complementary result allows us to calculate $\theta(G)$. These allow for $\chi(G)$ and $\theta(G)$ to be determined in polynomial time; see Corollary \ref{prop:chromaticnumber_polynomialtime}. 

\begin{theorem}
\label{prop:chromaticnumber}
    If $G \in $ \hng, then $\chi(G) = \omega(G)$ unless $G$ comprises an induced $C_5$ and an isolated set of vertices of types either (i) a subset of $\{c_1\}$, $\{c_1, c_3\},$ $\{c_1, c_4\}$, and $\{c_1, c_2, c_3, c_4, c_5\}$,  or (ii) both $\{c_1, c_2, c_4\}$ and $\{c_1, c_2, c_3, c_5\}$, up to symmetry. If either of these hold, $\chi(G) = \omega(G) + 1$. 
\end{theorem}
\begin{proof}
    Let $G \in$ \hngsp and suppose $G$ is imperfect. By Theorem \ref{prop:apex_perfect}, $G$ is apex-perfect, and contains at least one induced $C_5$ (note again that $G$ cannot contain any of the other forbidden induced subgraphs of perfect graphs). Let $C = \{c_1, c_2, c_3, c_4, c_5\}$ induce a copy of $C_5$ in $G$ with edge set $\{c_1c_2, c_2c_3, c_3c_4, c_4c_5, c_1c_5\}$ and let $D = V(G) - C$. Suppose $D$ matches none of the conditions given in the theorem and $|D| \ge 2$ (the result is straightforward to check if $|D| = 1$). If $\omega(D) = 1$ (i.e., $D$ induces a stable set), then by Proposition \ref{prop:allowable_vertex_list}, it follows that $D$ comprises vertices of types $\{c_1, c_2\}$ and $\{c_1, c_2, c_4\}$, up to symmetry. The argument in the following paragraph holds. Note that by the proof of Theorem \ref{prop:apex_perfect}, $D$ induces a perfect graph. 

    First, suppose $D$ has only vertices of types $\{c_1, c_2\}$ and $\{c_1, c_2, c_4\}$, up to symmetry. If so, a minimum proper coloring of $D$ can be extended to a proper coloring of $G$ using $\chi(D) + 2$ colors by giving vertices $c_1$ and $c_2$ new colors and coloring $c_4$ with $c_1$'s color, $c_5$ with $c_2$'s color, and $c_3$ with any color used in $D$. Hence $\chi(G) \le chi(D) + 2$. However, since $D$ is nonempty and complete to $c_1$ and $c_2$, it follows that $\omega(G) = \omega(D) + 2$. Since $D$ is perfect, we have $\chi(G) \ge omega(G) = \omega(D) + 2 = \chi(D) + 2$, and conclude that $\chi(G) = \omega(G)$. 

    Second, suppose $D$ may also have vertices of types $\{c_1, c_3\}$ and $\{c_1, c_2, c_3\}$, up to symmetry. If $D$ contains a vertex of type $\{c_1, c_3\}$, then $D$ contains only vertices of types $\{c_1, c_3\}, $ $\{c_1, c_4\},$ $\{c_1, c_3, c_4\},$ and $\{c_3, c_4\}$, of which there must be at least two vertices of types $\{c_1, c_3, c_4\},$ and $\{c_3, c_4\}$ to produce a clique number at least two. Vertices of types $\{c_1, c_3\}$ and $\{c_1, c_4\}$ can be freely colored with the same color as a vertex of type $\{c_1, c_3, c_4\}$ or $\{c_3, c_4\}$, so the above argument holds to show $\chi(G) = \omega(G)$. 
    
    If $D$ contains a vertex of type $\{c_1, c_2, c_3\}$, then $D$ comprises either (i) vertices of types $\{c_1, c_2\}$, $\{c_1, c_2, c_3\}$, $\{c_1, c_2, c_4\}$, $\{c_1, c_2, c_5\}$ or (ii) vertices of types $\{c_2, c_3\}$, $\{c_1, c_2, c_3\}$, $\{c_2, c_3, c_4\}$, $\{c_2, c_3, c_5\}$. Assuming the former, without loss of generality, $D$ is complete to $c_1$ and $c_2$, so a minimum proper coloring of $D$ extends to $G$ with two more colors: one is used for $c_1$ and $c_3$, one is used for $c_2$ and $c_4$, and $c_5$ can be colored with the same color as a vertex of type $\{c_1, c_2, c_3\}$. Note the coloring is proper since $D$ induces a clique, so no neighbor of $c_5$ in $D$ shares the same color as a vertex of type $\{c_1, c_2, c_3\}$. The above argument holds to show that $\chi(G) = \omega(G)$. 

    Third, suppose $D$ may also have vertices of types $\{c_1\}$ and $\{c_1, c_2, c_3, c_4\}$, up to symmetry. If $D$ contains a vertex of type $\{c_1\}$ not in a maximum clique, then $D$ without its vertices of type $\{c_1\}$ is as described in one of the above three paragraphs. Since the vertices of type $\{c_1\}$ are not adjacent to one another and not in a maximum clique in $D$, they can be added back into any proper coloring without using new colors. If a vertex of type $\{c_1\}$ is in a maximum clique, then since $\omega(D) \ge 2$, this maximum clique must also have one vertex of type $\{c_3, c_4\}$. Recall from the proof of Theorem \ref{prop:FIS_for_1_hng} that this maximum clique cannot contain two vertices of type $\{c_3, c_4\}$. Hence $\omega(D) = 2$, so all vertices of types $\{c_3, c_4\}$ and $\{c_1, c_3, c_4\}$ are isolated from one another. Here $G$ can be properly colored with three colors: the first is used on $c_2$, $c_5$, and all vertices of types $\{c_3, c_4\}$ and $\{c_1, c_3, c_4\}$; the second is used on $c_3$ and all vertices of type $\{c_1\}$; and the third is used on $c_1$ and $c_4$. 

    If $D$ contains a vertex of type $\{c_1, c_2, c_3, c_4\}$ and no vertex of type $\{c_2, c_3, c_5\}$, then by Proposition \ref{prop:allowable_vertex_list} it follows that $c_5$ has no neighbors in $D$. Hence any proper coloring of $D$ can be extended to a proper coloring of $G$ with two new colors, one for $c_1$ and $c_3$ and the other for $c_2$ and $c_4$, with a color from $D$ given to $c_5$. Furthermore, Proposition \ref{prop:allowable_vertex_list} implies that $\omega(G) = \omega(D) + 2$, so it follows that $\chi(G) = \omega(G)$, as reasoned above. If $D$ does have a vertex of type $\{c_2, c_3, c_5\}$, then by Proposition \ref{prop:allowable_vertex_list} all vertices in $D$ are adjacent to $c_2$ and $c_3$. As a result $\omega(G) = \omega(D) + 2$. Moreover, any minimum proper coloring of $D$ can be extended to a proper coloring of $G$ by giving $c_2$ and $c_4$ one new color and $c_1$ and $c_3$ a second new color. If there exists a vertex $v$ of type $\{c_1, c_2, c_3, c_4\}$ adjacent to all vertices of type $\{c_2, c_3, c_5\}$, then $v$ dominates $D$ and its color can be given to $c_5$ so that the coloring on $G$ is proper. Otherwise every vertex of type $\{c_1, c_2, c_3, c_4\}$ has a non-neighbor of type $\{c_2, c_3, c_5\}$. If there are two vertices of type $\{c_1, c_2, c_3, c_4\}, $ then $G$ contains $\overline{F_{25}}$ or $\overline{F_{26}}$ as an induced subgraph. If there are two adjacent vertices of type $\{c_2, c_3, c_5\}$, then these two vertices together with $c_1, c_2, c_5$, and the vertex of type $\{c_1, c_2, c_3, c_4\}$ induce an element of $\mathcal{F}$. Otherwise all vertices of type $\{c_2, c_3, c_5\}$ are non-adjacent to one another, and the vertex of type $\{c_1, c_2, c_3, c_4\}$ is adjacent to at least one of them. Here, $\omega(D) = 2, \omega(G) = 4, \chi(D) = 2$, and $\chi(G) = 4$.

    Lastly, suppose $D$ also has vertices of types $\{\emptyset\}$ and/or $\{c_1, c_2, c_3, c_4, c_5\}$. We showed in the proof of Theorem \ref{prop:apexperfect_cliquenumber} that a vertex of type $\{\emptyset\}$ is never in a maximum clique in $G$ or $D$, and a vertex of type $\{c_1, c_2, c_3, c_4, c_5\}$ is in all maximum cliques. Let $G'$ be the induced subgraph with no vertices of types $\{\emptyset\}$ or $\{c_1, c_2, c_3, c_4, c_5\}$, so by the above arguments, $\chi(G') = \omega(G')$. Suppose $D$ has $k$ vertices of type $\{c_1, c_2, c_3, c_4, c_5\}$. A minimum proper coloring of $G'$ can be extended to a proper coloring of $G$ with $\chi(G') + k$ colors, and $\omega(G) = \omega(G') + k$. Thus $\chi(G) \ge \omega(G) = \omega(G') + k = \chi(G') + k \ge \chi(G)$, so we conclude $\chi(G) = \omega(G)$. 

    For the second claim, let $G$ comprise an induced $C_5$ and a stable set $D$ of vertices of types either (i) a subset of $\{c_1\}$, $\{c_1, c_3\},$ $\{c_1, c_4\}$, and $\{c_1, c_2, c_3, c_4, c_5\}$,  or (ii) both $\{c_1, c_2, c_4\}$ and $\{c_1, c_2, c_3, c_5\}$, up to symmetry. In either case, if $D$ is nonempty, then $\omega(G) = 3$ and $\chi(G) = 4$. If $D$ is empty, then $G$ is isomorphic to $C_5$ and $\chi(G) = 3$. 
\end{proof}

\begin{corollary}
\label{prop:chromaticnumber_polynomialtime}
    If $G \in $ \hng, then $\chi(G)$ and $\theta(G)$ can be determined in polynomial time.
\end{corollary}
\begin{proof}
    By Corollary \ref{prop:cliquenumber_polynomialtime}, we can determine $\omega(G)$ in polynomial time. Subsequently, in $O(n^5)$ time, identify a subset of vertices $\{v_1, v_2, v_3, v_4, v_5\}$ inducing $C_5$, and check the condition in Theorem \ref{prop:chromaticnumber} to determine if $\chi(G) = \omega(G)$ or $\chi(G) = \omega(G) + 1$. 
\end{proof}


\begin{thebibliography}{99}

\frenchspacing

\bibitem{AoHa13} M. Aouchiche and P. Hansen, A survey of Nardhaus-Gaddum type relations, Discrete Appl. Math. \textbf{161} (4-5) (2013), 466--546.

\bibitem{BSZ} R. Behr, V. Sivaraman, and T. Zaslavsky,
\newblock Mock threshold graphs,
\newblock \emph{Discrete Math.} {\bf 341} (2018), 2159--2178. 

\bibitem{LWB} L. W. Beineke, 
\newblock Characterizations of derived graphs, 
\newblock \emph{J. Combin. Theory} {\bf 9} (1970), 129--135.

\bibitem{Bl93} Z. Bl\'{a}zsik et al.,
\newblock Graphs with no induced $C_4$ and $2K_2$, 
\newblock \emph{Discrete Math.} {\bf 115} (1993), 51--55. 

\bibitem{Be61} C. Berge, 
\newblock F\"{a}rburg von Graphen, deren s\"{a}mtlitche bzw. deren ungerade Kreise starr sind, 
\newblock \emph{Wiss. Z. Martin-Luther-Univ. Halle-Wittenberg Math.-Natur. Reihe} {\bf 10} (1961), 114 (in German). 

\bibitem{BLJ} A. Brandst\"adt, V. B. Le, and J. P. Spinrad, 
\newblock Graph Classes: A Survey, 
\newblock SIAM Monographs Discrete Math. Appl., Society for Industrial and Applied Mathematics, Philadelphia (1999). 

\bibitem{ChCoTr16}  C. Cheng, K. L. Collins, and A. N. Trenk,
\newblock Split graphs and Nordhaus--Gaddum graphs,
\newblock \emph{Discrete Math.} {\bf 339} (2016), 2345--2356.

\bibitem{Ch05} M. Chudnovsky et al., 
\newblock Recognizing Berge graphs, 
\newblock \emph{Combinatorica} {\bf 25} (2005), 143--186. 

\bibitem{CRST} M. Chudnovsky, N. Robertson, P. Seymour, and R. Thomas, 
\newblock The strong perfect graph theorem,
\newblock \emph{Ann. Math.} {\bf 164} (2006) (1), 51--229.

\bibitem{ChHa73} V. Chv\'{a}tal and P. L. Hammer, 
\newblock Set-packing and threshold graphs, 
\newblock Univ. Waterloo Res. Rep. CORR 73-21 (1973).

\bibitem{ChHa77} V. Chv\'atal and P.L. Hammer,
\newblock Aggregation of inequalities in integer programming,
\newblock in Annals of Discrete Mathematics 1: Studies in Integer Programming, ed. P. L. Hammer, E. I. Johnson, B. H. Korte, and G. L. Nemhauser, North Holland (1977), 145-162.

\bibitem{CoTr13}  K. L. Collins and A.N. Trenk,
\newblock Nordhaus-Gaddum Theorem for the Distinguishing Chromatic Number,
\newblock \emph{Electron. J. of Combin.} {\bf 20} (3) (2013). 

\bibitem{RD} R. Diestel, 
\newblock Graph Theory, 4th ed., 
\newblock in: Grad. Texts in Math. {\bf 173}, Springer, Heidelberg (2010).

\bibitem{Di61} G. A. Dirac,
\newblock On rigid circuit graphs, 
\newblock \emph{Abh. Math. Sem. Univ. Hamburg} {\bf 25} (1961), 71--76.

\bibitem{CFGS} R. Faudree, E. Flandrin, and Z. Ryj\'a\v{c}ek, 
\newblock Claw-free graphs--a survey, 
\newblock \emph{Discrete Math.} {\bf 164} (1997), 87--147.

\bibitem{Fi66} H. J. Finck, 
\newblock On the chromatic number of a graph and its complements,
\newblock Theory of Graphs, Proceedings of the Colloquium, Tihany, Hungary (1966), 99--113. 

\bibitem{FoHa77}  S. F\"{o}ldes and P. Hammer, 
\newblock Split graphs,
\newblock \emph{Congr. Numer.} {\bf 19} (1977), 311--315.

\bibitem{MCG} M. C. Golumbic, 
\newblock Algorithmic Graph Theory and Perfect Graphs,
\newblock Academic Press, New York, 1980. 2nd ed., Ann. Discrete Math. {\bf{57}} (2004).

\bibitem{GrLoSc84} M. Gr\"{o}tschel, L. Lov\'{a}sz, and A. Schrijver, 
\newblock Polynomial algorithms for perfect graphs, 
\newblock Ann. Discrete Math. {\bf{21}} (1984), 325-356. 

\bibitem{HaSi81}  P. Hammer and B. Simeone,
\newblock The splittance of a graph,
\newblock \emph{Combinatorica} {\bf 1} (1981), 275--284.

\bibitem{Ha85} R. B. Hayward, 
\newblock Weakly triangulated graphs, 
\newblock \emph{J. Combin. Theory Ser. B} {\bf 39} (3) (1985), 200--208.

\bibitem{LiPoSi19} B. Litjens, S. Polak, and V. Sivaraman, 
\newblock Sum-perfect graphs, 
\newblock \emph{Discret. Appl. Math.} {\bf 259} (2019), 232--239. 

\bibitem{LL} L. Lov\'{a}sz, 
\newblock A characterization of perfect graphs, 
\newblock \emph{J. Combin. Theory Ser. B} \textbf{13} (1972), 95--98.

\bibitem{MaPe95} N.V.R. Mahadev and U.N. Peled,
\newblock \emph{Threshold Graphs and Related Topics,}
\newblock in: Ann. Discrete Math. {\bf 56}, North-Holland, Amsterdam (1995).

\bibitem{NoGa56} E. A. Nordhaus and J. W. Gaddum, 
\newblock On complementary graphs, 
\newblock \emph{Amer. Math. Monthly} {\bf 63} (1956), 175--177. 

\bibitem{ScSe20} A. Scott and P. Seymour, 
\newblock A Survey of $\chi$-boundedness, 
\newblock \emph{J. Graph Theory} {\bf 95} (3) (2020), 473--504. 

\bibitem{StTu08} C. L. Starr and G. E. Turner III, 
\newblock Complementary graphs and the chromatic number, 
\newblock \emph{Missouri J. Math Sci.} {\bf 20} (1) (2008), 19--26. 

\bibitem{DBW}  D. B. West, 
\newblock Introduction to Graph Theory, 2nd ed., 
\newblock Prentice Hall, Upper Saddle River, N.J. (2001).

\bibitem{YCC}  J.-H. Yan, J.-J. Chen, and G. J. Chang,
\newblock Quasi-threshold graphs, 
\newblock \emph{Discrete Appl. Math.} {\bf 69} (3) (1996), 247--255.

\end{thebibliography}
\end{document}